\documentclass[journal]{IEEEtran}

\usepackage{etex}
\usepackage[vlined,ruled]{algorithm2e}
\usepackage[dvipsnames]{xcolor}
\usepackage{pgfpages}
\usepackage[cmex10]{amsmath}		
\usepackage{amssymb, mathrsfs}
\usepackage{cite}
\usepackage{url}

\usepackage{dsfont}
\usepackage{siunitx}
\usepackage{verbatim}									
\usepackage{mathtools}									
\usepackage{siunitx}									
\usepackage{mdframed}

\usepackage{graphicx}
\usepackage{epstopdf}

\usepackage[font=small]{caption}
\usepackage{subcaption}
\usepackage[activate={true,nocompatibility},final,tracking=true,kerning=true,spacing=true,factor=1100,stretch=10,shrink=10]{microtype}

\usepackage{array}
\usepackage{multirow}
\usepackage{enumerate}
\usepackage{booktabs}

\graphicspath{{./fig/}}

\usepackage{pgfplots}
\pgfplotsset{width=8cm,compat=newest}
\newcommand*{\damping}{0.006}%
\newcommand*{\freq}{25}%
\pgfmathsetmacro{\freqd}{sqrt(1-(\damping)^2)*\freq}%

\pgfplotsset{
    standard/.style={
    axis x line=middle,
    axis y line=middle,
    enlarge x limits=0.15,
	enlarge y limits=0.15,
	every axis plot post/.style={mark options={fill=black}},
	}
}

\pgfplotsset{%
    ,compat=1.12
    ,every axis x label/.style={at={(current axis.right of origin)},anchor=north west}
    ,every axis y label/.style={at={(current axis.above origin)},anchor=north east}
    }

\usepackage{tikz}
\usetikzlibrary{arrows}
\usetikzlibrary{decorations.pathmorphing}
\usetikzlibrary{decorations.markings}
\usetikzlibrary{snakes}
\usetikzlibrary{matrix}
\usetikzlibrary{calc}
\usetikzlibrary{patterns}
\usetikzlibrary{positioning}
\usetikzlibrary{shapes}
\usetikzlibrary{shapes.symbols}
\usetikzlibrary{shapes.geometric}
\usetikzlibrary{shadows.blur}
\usetikzlibrary{shapes.arrows}
\usetikzlibrary{shapes.multipart}
\usetikzlibrary{shapes.callouts}
\usetikzlibrary{shapes.misc}

\tikzstyle{every node}=[font=\small]
\tikzstyle{every path}=[line width=0.8pt,line cap=round,line join=round]

\usepackage[american,smartlabels]{circuitikz}
\ctikzset{bipoles/thickness=1}
\ctikzset{bipoles/length=0.8cm}
\ctikzset{bipoles/diode/height=.375}
\ctikzset{bipoles/diode/width=.3}
\ctikzset{tripoles/thyristor/height=.8}
\ctikzset{tripoles/thyristor/width=1}
\ctikzset{bipoles/vsourceam/height/.initial=.7}
\ctikzset{bipoles/vsourceam/width/.initial=.7}

\newcommand{\real}{\mathbb{R}}

\newcommand{\setdef}[2]{\{#1 \;|\; #2\}}

\DeclareMathOperator*{\minimize}{minimize} 									
\DeclareMathOperator{\subto}{subject~to}

\newcommand{\vect}[1]{\mathbbold{#1}}

\newcommand{\vzeros}[1][]{\vect{0}_{#1}}
\DeclareSymbolFont{bbold}{U}{bbold}{m}{n}
\DeclareSymbolFontAlphabet{\mathbbold}{bbold}

\newcommand{\map}[3]{#1: #2 \rightarrow #3}

\newcommand{\tr}{\color{red}}

\newcommand{\tb}{\color{blue}}

\newcommand{\T}{\mathsf{T}} 

\newcommand\oprocendsymbol{\hbox{$\square$}}
\newcommand\oprocend{\relax\ifmmode\else\unskip\hfill\fi\oprocendsymbol}

\newtheorem{theorem}{Theorem}[section]

\newtheorem{remark}[theorem]{Remark}
\newtheorem{assumption}[theorem]{Assumption}
\newtheorem{example}[theorem]{Example}

\newenvironment{pfof}[1]{\vspace{1ex}\noindent{\itshape Proof of
    #1:}\hspace{0.5em}} {\hfill\oprocend\vspace{1ex}}

\makeatletter
\newcommand{\vast}{\bBigg@{4}}
\newcommand{\Vast}{\bBigg@{5}}

\newcommand{\define}{\ensuremath{\triangleq}}

\usepackage{enumitem}

\usepackage{enumerate}
\usepackage{cuted}
\usepackage{arydshln}

\renewcommand{\tb}{\color{black}}

\begin{document}
\title{{\tb Analysis and Synthesis of Low-Gain Integral Controllers for Nonlinear Systems}}

\author{John W. Simpson-Porco,~\IEEEmembership{Member,~IEEE}\\
\thanks{J.~W.~Simpson-Porco is with the Department of Electrical and Computer Engineering, University of Toronto, 10 King's College Road,
Toronto, ON, M5S 3G4, Canada. Email: {\tt jwsimpson@ece.utoronto.ca}.}
\thanks{{\tb This work was supported by NSERC Discovery Grant RGPIN-2017-04008.}}
}

\markboth{Submitted for Publication. This version: \today}%
{Submitted to IEEE Transactions on Automatic Control. This version: \today}

\maketitle

\begin{abstract}
Relaxed conditions are given for stability of a feedback system consisting of an exponentially stable multi-input multi-output nonlinear plant and an integral controller. Roughly speaking, it is shown that if the composition of the plant equilibrium input-output map and the integral feedback gain is infinitesimally contracting, then the closed-loop system is exponentially stable if the integral gain is sufficiently low. {\tb The main result is illustrated with an application arising in frequency control of AC power systems. We demonstrate how the contraction condition can be checked computationally via semidefinite programming, and how integral gain matrices can be synthesized via convex optimization to achieve robust $\mathscr{L}_2$ performance in the presence of nonlinearity and uncertainty.} 
\end{abstract}

\begin{IEEEkeywords}
Nonlinear control, integral control, output regulation, linear matrix inequalities (LMIs), slow integrators, singular perturbation
\end{IEEEkeywords}

\section{Introduction}
\label{Sec:Introduction}


Tracking and disturbance rejection in the presence of model uncertainty is one of the fundamental purposes of automatic control. {\tb A case which commonly occurs in engineering practice is that the system one wishes to control is complex, and no accurate dynamic model is available, but it is however known that the system is stable and is responsive to control inputs. This stability may be inherent to the system, or may have been achieved through a preliminary stabilizing feedback design. A modest and practical design goal is then simply to improve reference tracking and disturbance rejection performance via a supplementary integral controller, without compromising system stability.
One concrete example of this problem is that of frequency control for large-scale AC power systems, where the Automatic Generation Control (AGC) system asymptotically rebalances load and generation via integral action. In practical AGC systems, the integral gain is set very low, to ensure that the uncertain power system dynamics are not destabilized by the supplementary feedback \cite{JWSP:20d}.}

{\tb In the SISO LTI case, tuning such a supplementary integral loop requires only that one know the sign of the plant DC gain}; the general MIMO LTI case is slightly more subtle. Consider the continuous-time LTI state-space model
\begin{equation}\label{Eq:LTI}
\begin{aligned}
\dot{x} &= Ax + Bu + B_ww\\
e &= Cx + Du + D_ww
\end{aligned}
\end{equation}
with state $x \in \real^n$, control input $u \in \real^m$, constant disturbance/reference signal $w \in \real^{n_w}$, and error output $e \in \real^p$. {\tb For the reasons described previously, assume that $A$ is Hurwitz.} One interconnects the system \eqref{Eq:LTI} with the integral controller
\begin{equation}\label{Eq:LTIIntegral}
\begin{aligned}
\dot{\eta} &= -\varepsilon e, \qquad u = K\eta,
\end{aligned}
\end{equation}
where $K \in \real^{m \times p}$ is a gain matrix and $\varepsilon > 0$. Let 
\begin{equation}\label{Eq:DCGains}
G(0) = -CA^{-1}B+D, \quad G_w(0) = -CA^{-1}B_w + D_w
\end{equation}
denote the DC gain matrices of \eqref{Eq:LTI} from $u$ and $w$ to $e$. Implicit in the proof of \cite[Lemma 3]{EJD:76} is the following: if there exists $K$ such that $-G(0)K$ is Hurwitz, then the feedback system \eqref{Eq:LTI}--\eqref{Eq:LTIIntegral} is internally exponentially stable for {\tb sufficiently} small $\varepsilon > 0$. {\tb The existence of $K$ such that $-G(0)K$ is Hurwitz is equivalent to} $G(0)$ having full row rank, and indeed a suitable gain design is $K = G(0)^{\dagger}$, as used in \cite[Lemma 3]{EJD:76}. This result was stated succinctly in \cite[Theorem 3]{MM:85}; see also \cite[Lemma 1, A.2, A.3]{DEM-EJD:89} for further details.

{\tb From a singular perturbation point of view \cite{PVK-HKK-JO:99}, the low integral gain $\varepsilon$ induces a \emph{time-scale separation} in the system \eqref{Eq:LTI}--\eqref{Eq:LTIIntegral}, and this is the perspective we will exploit going forward. The key intuition is as follows. If $\varepsilon$ is sufficiently small in \eqref{Eq:LTIIntegral}, then $\eta(t)$ and $u(t)$ change very slowly. Relatively speaking then, the plant dynamics \eqref{Eq:LTI} are fast, and the output signal $e(t)$ will be well approximated by the \emph{quasi steady-state} relationship $e(t) = G(0)u(t) + G_w(0)w$. Substituting this into \eqref{Eq:LTIIntegral} leads to the simplified \emph{slow time-scale} dynamics
\begin{equation}\label{Eq:SlowDynamics}
\dot{\eta} = - \varepsilon G(0)K\eta - \varepsilon G_w(0)w.
\end{equation}
In other words, the model \eqref{Eq:SlowDynamics} describes the closed-loop dynamics when one ignores the fast plant dynamics. By inspection, \eqref{Eq:SlowDynamics} is internally exponentially stable if and only if $-G(0)K$ is Hurwitz, which is again the Davison/Morari result.}


Extensions of this LTI result to Lur'e-type systems \cite{TF-HL-EPR:03,CG-HL-ST:17}, {\tb discrete-time systems \cite{HL-ST:97}, }and to distributed-parameter systems \cite{ATJ-VA-VDSM-CX:19} have been pursued. In the full nonlinear setting, the most well-known result is due to Desoer and Lin \cite{CD-CL:85}, who proved that if the \emph{equilibrium input-to-error map} of the plant is strongly monotone, then a similar low-gain stability result holds; a related condition was recently also used in \cite[Equation (21)]{XH-HKH-YS:19}. When specialized for LTI systems, the Desoer--Lin condition states that $G(0) + G(0)^{\T}$ should be positive definite; it therefore does not properly generalize the Davison/Morari result. It appears the only attempt to close this gap was reported in \cite{HCT:92}, where Tseng proposed a design based on inverting the Jacobian of the plant equilibrium input-to-error map. This recovers Davison's special design $K = G(0)^{-1}$ in the (square) LTI case, but in general yields a very complicated nonlinear feedback. In sum, the available {\tb low-gain integral control} results in the literature {\tb for nonlinear systems} do not reduce as expected in the LTI case, and the literature lacks systematic procedures for constructing low-gain integral controllers for nonlinear and uncertain systems. One goal of this paper is to close this gap.

While the design of traditional low-gain tracking controllers is important in and of itself, another source of recent interest in such low-gain methods in a nonlinear context has arisen from the study of feedback-based optimizing controllers for dynamic systems; see \cite{LSPL-ZEN-EM-JWSP:18e,LSPL-JWSP-EM:18l, MC-EDA-AB:18, SM-AH-SB-GH-FD:18, MC-JWSP-AB:19c} for various recent works. In this line of work, the controller does not attempt to track an explicit reference, but instead attempts to drive the system towards an optimal equilibrium point in the presence of an unmeasured exogenous disturbance. {\tb As such controllers are based on regulating a suitable measure of sub-optimality to zero}, the results we develop will also be applicable to this class of ``tracking-adjacent'' problems.

\medskip

\paragraph*{Contributions} 

The broad goals of this paper are (i) to understand when low-gain MIMO integral feedback can be applied to a MIMO nonlinear system, and (ii) to leverage modern robust control tools for the analysis and design of such supplementary loops. These goals are largely inspired by practical problems in power system control, and by the foundational paper \cite{EJD:76}, which provided a constructive solution in the LTI case. As a result of these goals, this work is somewhat disjoint from the modern literature on output regulation (see \cite{DA-LP:17,XH-HKH-YS:19,MB-LM:19a,MB-LM:19b} for recent contributions) where the focus is on quite different issues, such as nonlinear stabilization, practical vs. asymptotic regulation, and the construction of internal models. This work is therefore best understood as a continuation of the line of research in \cite{EJD:76,CD-CL:85,HCT:92}.

There are two main contributions. First, in Section \ref{Sec:Desoer} we present a generalization the main result of \cite{CD-CL:85}, providing relaxed conditions on the plant's equilibrium input-to-error map which ensure closed-loop stability under low-gain integral control. The main idea is to impose that the reduced time-scale dynamics be infinitesimally contracting, which ensures the existence of a unique and exponentially stable equilibrium point for all  constant exogenous disturbances. Unlike the conditions reported in \cite{CD-CL:85,HCT:92}, this condition recovers the Davison/Morari result that $-G(0)K$ should be Hurwitz when restricted to the LTI case, and allows for additional flexibility over \cite{CD-CL:85,HCT:92}. We apply the results to show stability of a nonlinear frequency regulation scheme for AC power systems. Second, in Section \ref{Sec:SDP} we describe how semidefinite programming can be used for certification of stability under low-gain integral control, as well as for direct convex synthesis of integral gain matrices which achieve robust performance. These results apply equally in the nonlinear or in the uncertain linear contexts, and are illustrated via academic examples. 


\subsection{Notation} 

For column vectors $x_1,\ldots x_n$, $\mathrm{col}(x_1,\ldots,x_n)$ denotes the associated concatenated column vector. The identity matrix of size $n$ is $I_n$, and $\vzeros[n]$ is the zero vector of dimension $n$.
%
%
The notation $P \succ 0$ (resp. $P \succeq 0$) means that the matrix $P$ is symmetric and positive (semi)definite. In any expression of the form $(\star)^{\T}XY$ or $\left[\begin{smallmatrix}X & Y \\
(\star)^{\T} & Z\end{smallmatrix}\right]$, $(\star)$ is simply an abbreviation for $Y$.
For matrices $X_1,\ldots,X_n$, $\mathrm{diag}(X_1,\ldots,X_n)$ denotes the associated block diagonal matrix. Given two $2\times 2$ block partitioned matrices $X = \left[\begin{smallmatrix}X_{11} & X_{12} \\ X_{21} & X_{22}\end{smallmatrix}\right]$ and $Y = \left[\begin{smallmatrix}Y_{11} & Y_{12} \\ Y_{21} & Y_{22}\end{smallmatrix}\right]$, we define the diagonal augmented matrix
\[
\mathrm{daug}(X,Y) = {\small \left[\begin{array}{@{}cc|cc@{}}
X_{11} & 0 & X_{12} & 0\\
0 & Y_{11} & 0 & Y_{12}\\
\hline
X_{21} & 0 & X_{22} & 0\\
0 & Y_{21} & 0 & Y_{22}
\end{array}\right]}.
\]
Given a scalar-valued function $V(x,y)$, $\nabla_x V(x,y)$ and $\nabla V_{y}(x,y)$ denote its gradients with respect to $x$ and $y$, respectively. The space $\mathscr{L}_{2}^{p}[0,\infty)$ denotes the set of measurable maps $\map{f}{\real}{\real^p}$ which are zero for $t < 0$ with $\tau \mapsto \|f(\tau)\|_{2}^2$ being integrable over $[0,\infty)$, and $\mathscr{L}_{2e}^{p}[0,\infty)$ denotes the associated extended signal space where $\tau \mapsto \|f(\tau)\|_{2}^2$ is integrable over $[0,T]$ for all $T > 0$; see, e.g., \cite{AJvdS:16} for details. For $f \in \mathscr{L}_{2}^{p}[0,\infty)$, we let $\|f\|_{\mathscr{L}_2} = (\int_{0}^{\infty}\|f(\tau)\|_2^2 \,\mathrm{d}\tau)^{1/2}$.

%


\section{Problem Setup and Assumptions}
\label{Sec:ProblemSetup}

We consider a physical plant which is described by a finite-dimensional nonlinear time-invariant state-space model
\begin{equation}\label{Eq:Plant}
\begin{aligned}
\dot{x}(t) &= f(x(t),u(t),w), \quad x(0) = x_0\\
e(t) &= h(x(t),u(t),w)
\end{aligned}
\end{equation}
where $x(t) \in \real^n$ is the state with initial condition $x_0$, $u(t) \in \real^m$ is the control input, $e(t) \in \real^p$ is the error to be regulated to zero, and $w \in \real^{n_w}$ is a vector of \emph{constant} reference signals, disturbances, and unknown parameters.\footnote{Asymptotically constant references/disturbances are treated similarly \cite{KHH:00}.} The maps $f$ and $h$ are defined on a domain of interest $D_{x} \times D_{u} \times D_{w} \subseteq \real^{n} \times \real^{m} \times \real^{n_w}$. For fixed $w$, the possible equilibrium state-input-error triplets $(\bar{x},\bar{u},\bar{e})$ are determined by the algebraic equations
\[
\begin{aligned}
\vzeros[n] &= f(\bar{x},\bar{u},w), \qquad \bar{e} = h(\bar{x},\bar{u},w).
\end{aligned}
\]
{\tb
We next lay out our main assumptions on the plant \eqref{Eq:Plant}.

\smallskip

\begin{assumption}[\bf Plant Assumptions]\label{Ass:Plant}
For \eqref{Eq:Plant} there exist sets $\mathcal{X} \subseteq D_x$ and $\mathcal{I} \subseteq D_u \times D_w$ such that
\begin{enumerate}[label=(A\arabic*)]
\item \label{Ass:Plant-0} $f$, $h$, and all associated Jacobian matrices are Lipschitz continuous on $\mathcal{X}$, uniformly with respect to $(u,w) \in \mathcal{I}$;
\item \label{Ass:Plant-1} there exists a continuously differentiable map $\map{\pi_{x}}{\mathcal{I}}{\mathcal{X}}$ which is Lipschitz continuous on $\mathcal{I}$ and satisfies
\[
\vzeros[n] = f(\pi_{x}(u,w),u,w), \qquad \text{for all}\,\,\, (u,w) \in \mathcal{I}.
\]
\item \label{Ass:Plant-3} the equilibrium $\bar{x} = \pi_x(u,w)$ is exponentially stable, uniformly in $(u,w) \in \mathcal{I}$.
%
%
%
%
%
%
%
%

%
\end{enumerate}
\end{assumption}

Assumption \ref{Ass:Plant} essentially states that each constant input-disturbance pair $(u,w) \in \mathcal{I}$ yields a unique (at least locally on the set $\mathcal{X}$) exponentially stable equilibrium state $\bar{x} = \pi_x(u,w)$.
More specifically, \ref{Ass:Plant-1} states that the plant possesses an equilibrium which changes smoothly over the set of considered inputs, while \ref{Ass:Plant-0} and \ref{Ass:Plant-3} ensure sufficient model smoothness and uniform exponential stability of the equilibrium. The specific conditions in \ref{Ass:Plant-0} and \ref{Ass:Plant-3} have been chosen to satisfy the conditions of a relatively standard converse Lyapunov theorem \cite[Lemma 9.8]{HKK:02}; variations on \ref{Ass:Plant-0} and \ref{Ass:Plant-3} are therefore likely possible.} When restricted to the LTI case in \eqref{Eq:LTI}, Assumption \ref{Ass:Plant} simply reduces to the matrix $A$ being Hurwitz, and we may select $\mathcal{X} = \real^n$ and $\mathcal{I} = \real^{m} \times \real^{n_w}$.

We call the map $\map{\pi}{\mathcal{I}}{\real^p}$ given by
\begin{equation}\label{Eq:DefofPi}
\pi(u,w) \define h(\pi_{x}(u,w),\bar{u},w)
\end{equation}
the \emph{equilibrium input-to-error map}. 
%
%
%
%
%
%
For notational use, we let
\[
\begin{aligned}
\mathcal{W} &= \setdef{w \in \real^{n_w}}{\text{there exists}\,\,u \in \real^{m}\,\,\text{s.t.}\,\,(u,w) \in \mathcal{I}}
\end{aligned}
\]
and for a given $w \in \mathcal{W}$ we let
\begin{equation}\label{Eq:Uw}
\mathcal{U}_{w} = \setdef{u \in \real^{m}}{(u,w) \in \mathcal{I}} \neq \emptyset
\end{equation}
denote the set of constant controls for which the equilibrium map $u \mapsto \pi_{x}(u,w)$ is defined.

\medskip

\begin{example}[\bf Illustration of Assumption \ref{Ass:Plant}] Consider the scalar dynamic system 
\[
\begin{aligned}
\dot{x} &= -\beta\sin(x) + u - w, \qquad e = h(x,u,w) = x,
\end{aligned}
\]
where $\beta > 0$. 
For fixed $\gamma \in [0,\tfrac{\pi}{2})$, define
\[
\begin{aligned}
\mathcal{I}(\gamma) &= \setdef{(u,w) \in \real^2}{-\beta\sin(\gamma) \leq  u-w \leq \beta\sin(\gamma)}\\
\mathcal{X}(\gamma) &= \setdef{x \in \real}{-\gamma \leq x \leq \gamma}.
\end{aligned}
\]
{\tb 
It follows then that one equilibrium is given by $\bar{x} = \pi_x(u,w) = \pi(u,w) = \arcsin((u+w)/\beta)$, and indeed $\pi_x$ is a continuously differentiable and Lipschitz continuous map from $\mathcal{I}(\gamma)$ into $\mathcal{X}(\gamma)$. Exponential stability of $\bar{x}$ can be certified with a simple quadratic Lyapunov function $V(x,u,w) = (x-\pi(u,w))^2$, with uniformity of stability in $(u,w)$ following by compactness of $\mathcal{I}(\gamma)$; the details are omitted. \hfill \oprocend
}
\end{example}

\medskip

\begin{figure}[ht!]
\begin{center}
\tikzstyle{block} = [draw, fill=white, rectangle, 
    minimum height=3em, minimum width=6em, blur shadow={shadow blur steps=5}]
    \tikzstyle{hold} = [draw, fill=white, rectangle, 
    minimum height=3em, minimum width=3em, blur shadow={shadow blur steps=5}]
    \tikzstyle{widehold} = [draw, fill=white, rectangle, 
    minimum height=3em, minimum width=5em, blur shadow={shadow blur steps=5}]
    \tikzstyle{dzblock} = [draw, fill=white, rectangle, minimum height=3em, minimum width=4em, blur shadow={shadow blur steps=5},
path picture = {
\draw[thin, black] ([yshift=-0.1cm]path picture bounding box.north) -- ([yshift=0.1cm]path picture bounding box.south);
\draw[thin, black] ([xshift=-0.1cm]path picture bounding box.east) -- ([xshift=0.1cm]path picture bounding box.west);
\draw[very thick, black] ([xshift=-0.5cm]path picture bounding box.east) -- ([xshift=0.5cm]path picture bounding box.west);
\draw[very thick, black] ([xshift=-0.5cm]path picture bounding box.east) -- ([xshift=-0.1cm, yshift=+0.4cm]path picture bounding box.east);
\draw[very thick, black] ([xshift=+0.5cm]path picture bounding box.west) -- ([xshift=+0.1cm, yshift=-0.4cm]path picture bounding box.west);
}]
\tikzstyle{sum} = [draw, fill=white, circle, node distance=1cm, blur shadow={shadow blur steps=8}]
\tikzstyle{input} = [coordinate]
\tikzstyle{output} = [coordinate]
\tikzstyle{pinstyle} = [pin edge={to-,thin,black}]

\begin{tikzpicture}[auto, scale = 0.6, node distance=2cm,>=latex', every node/.style={scale=1}]
    \node [block] (system) {{$
    \begin{aligned}
    \dot{x} &= f(x,u,w)\\
    e &= h(x,u,w)
    \end{aligned}$}};
    \node [widehold, below of=system, node distance=2cm, xshift=0.75cm] (controller)  {$\dot{\eta} = -\varepsilon e$};
    \node [hold, left of=controller, node distance=2.2cm] (nonlinear)  {$k(\cdot)$};
    \node [output, right of=system, node distance=3cm] (output) {};

    \coordinate [left of=system, node distance=2.5cm] (tmp);
    \draw [thick, -latex] (system) -- node [name=y] {$e$}(output);
    \draw [thick, -latex] (y) |- (controller);
    \draw [thick, -latex] (controller) -- node[above] {$\eta$} (nonlinear);
    \draw [thick, -latex] (nonlinear) -| (tmp) -- node[] {$u$} (system);
    
        \node [input, above of=system, node distance=1.2cm] (dist) {};
        \draw [thick,-latex] (dist) -- node[] {$w$} (system.north);
\end{tikzpicture}
\end{center}
\caption{{\tb Block diagram of plant with low-gain integral controller.}}
\label{Fig:LowGain}
\end{figure}
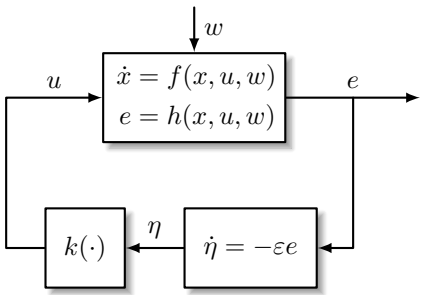

We are interested in the application of a pure integral feedback control scheme to \eqref{Eq:Plant} which acts on the error $e$ as
\begin{equation}\label{Eq:Integral}
\begin{aligned}
\dot{\eta} &= -\varepsilon e, \quad \eta(0) \in \real^{p}\\
u &= k(\eta)
\end{aligned}
\end{equation}
where $\map{k}{\real^p}{\real^m}$ is a feedback and $\varepsilon > 0$ is to be determined. We assume that
\smallskip
\begin{enumerate}[label=(A\arabic*)]\setcounter{enumi}{3}
\item \label{Ass:Controller-1} $k(\eta)$ is continuously differentiable and $L_k$-Lipschitz continuous on $\real^p$.
\end{enumerate}
\smallskip
The closed-loop system is the interconnection of the plant \eqref{Eq:Plant} and the controller \eqref{Eq:Integral}{\tb , and is shown in Figure \ref{Fig:LowGain}. }

Our goal is to give conditions under which \eqref{Eq:Plant},\eqref{Eq:Integral} possesses an exponentially stable equilibrium point for sufficiently small $\varepsilon > 0$ {\tb and for all disturbances $w \in \mathcal{W}$.}

\section{A Generalization of the Desoer-Lin Result on Low-Gain Integral Control}
\label{Sec:Desoer}

The main result of this section provides a generalization of the result of \cite{CD-CL:85}, where the monotonicity requirement on the equilibrium input-to-error map $u \mapsto \pi(u,w)$ is weakened to infinitesimal contraction \cite{WL-JJES:98, EDS:10, JWSP-FB:12za, FF-RS:14} of the vector field $\eta \mapsto -\pi(k(\eta),w)$. {\tb There are several motivations for working with contraction-type stability criteria in this context. First, contraction allows us to obtain stability results which are independent of the operating point, and hence independent of the exogenous disturbance $w$. Second, contractive systems possess simple Lyapunov functions, a fact we will exploit in the proof of Theorem \ref{Thm:Desoer}. Third, note that the vector field $\pi(k(\eta),w)$ has dimension equal to the number of regulated outputs; this dimension will be low in many practical problems. Contraction analysis can often be performed analytically for low-dimensional vector fields, particularly in the global setting (see the example in Section \ref{Sec:PowerSystemFrequencyControl}). Fourth and finally, contraction analysis is compatible with LMI-based analysis and design techniques \cite[Chap. 4.3]{SB-LEG-EF-VB:94}. \cite[Chap. 5]{AP-NVDW-HN:05}, which we will exploit in Section \ref{Sec:SDP}. These motivations aside, alternatives to contraction analysis are certainly possible; see Remark \ref{Rem:Alternatives} at the end of this section.
}


While multiple approaches to contraction analysis can be found in the literature, of varying sophistication, we will make use of the formulation based on the \emph{matrix measure}; this has proved sufficient for our applications of interest. Let $\|\cdot\|$ denote any vector norm on $\real^n$, with $\|\cdot\|$ also denoting the associated induced matrix norm. The \emph{matrix measure} associated with $\|\cdot\|$ is the mapping $\map{\mu}{\real^{n\times n}}{\real}$ defined by \cite[Chap. 2.2.2]{MV:02}
\[
\mu(A) \define \lim_{h\rightarrow 0^+}\tfrac{1}{h}(\|I_n+hA\|-1).
\]
Matrix measures associated with standard vector norms $\|\cdot\|_{1}$, $\|\cdot\|_{2}$, $\|\cdot\|_{\infty}$ (and their weighted variants) are all explicitly computable and have found substantial use in applications; see \cite{EDS:10,ZA-EDS:14} for clear summaries.

Infinitesimal contraction of a vector field is characterized via the matrix measure of its Jacobian matrix \cite{EDS:10}. Let $\mathcal{P}$ be a non-empty parameter set, and consider the dynamics
\begin{equation}\label{Eq:ParameterVectorField}
\dot{x}(t) = F(x(t),\alpha), \qquad x(0) = x_0, \quad \alpha \in \mathcal{P},
\end{equation}
where $\map{F}{\real^n \times \mathcal{P}}{\real^n}$ is continuously differentiable in its first argument. Let $\|\cdot\|$  be any vector norm. For a given $\alpha \in \mathcal{P}$, the system \eqref{Eq:ParameterVectorField} is \emph{infinitesimally contracting} with respect to $\|\cdot\|$ on a set $\mathcal{C}_{\alpha} \subseteq \real^n$ if there exists $\rho_{\alpha} > 0$ such that 
\begin{equation}\label{Eq:MatrixMeasureBound}
\mu(\tfrac{\partial F}{\partial x}(x,\alpha)) \leq -\rho_{\alpha}, \quad \text{for all}\,\, x \in \mathcal{C}_{\alpha}.
\end{equation}
If $\mathcal{C}_{\alpha}$ is a convex and forward-invariant set for \eqref{Eq:ParameterVectorField}, then \eqref{Eq:MatrixMeasureBound} guarantees that  \eqref{Eq:ParameterVectorField} possesses a unique equilibrium point $x^{\star} \in \mathcal{C}_{\alpha}$ towards which all trajectories with initial conditions $x_0 \in \mathcal{C}_{\alpha}$ will converge exponentially at rate $e^{-\rho_{\alpha}t}$ \cite[Thm. 1/2]{EDS:10}. 

As the parameter $\alpha$ varies over $\mathcal{P}$, it is generally the case that both the foward-invariant set $\mathcal{C}_{\alpha}$ and the contraction rate $\rho_{\alpha}$ will need to vary. To ensure a uniform rate, we will say that \eqref{Eq:ParameterVectorField} is \emph{uniformly infinitesimally contracting} with respect to $\|\cdot\|$ on a family of sets $\{\mathcal{C}_{\alpha}\}_{\alpha \in \mathcal{P}}$ if there exists $\rho > 0$ such that for all $\alpha \in \mathcal{P}$
\begin{equation}\label{Eq:MatrixMeasureBoundRobust}
\mu(\tfrac{\partial F}{\partial x}(x,\alpha)) \leq -\rho, \quad \text{for all}\,\, x \in \mathcal{C}_{\alpha}.
\end{equation}
We are now ready to state our main result.

\medskip

\begin{theorem}[\bf Relaxed Conditions for Exponential Stability under Low-Gain Integral Control]\label{Thm:Desoer}
Consider the plant \eqref{Eq:Plant} interconnected with the integral controller \eqref{Eq:Integral} under assumptions \ref{Ass:Plant-0}--\ref{Ass:Controller-1}. Define the reduced dynamics
\begin{equation}\label{Eq:NonlinearSlow}
\dot{\eta} = F_{\rm s}(\eta,w) \define -\pi(k(\eta),w), \quad \eta(0) =  \eta_0,
\end{equation}
and assume that 
\begin{enumerate}[label=(A\arabic*)]\setcounter{enumi}{4}
\item \label{Thm:Desoer-1} for each $w \in \mathcal{W}$ there exists a convex forward-invariant set $\mathcal{C}_{w}$ for \eqref{Eq:NonlinearSlow} such that $k(\mathcal{C}_{w}) \subseteq \mathcal{U}_{w}$. 
\item \label{Thm:Desoer-2} the system \eqref{Eq:NonlinearSlow} is uniformly infinitesimally contracting with respect to some norm $\|\cdot\|_{\rm s}$ on $\{\mathcal{C}_{w}\}_{w \in \mathcal{W}}$. 
\end{enumerate}
Then there exists $\varepsilon^{\star} > 0$ such that for any $\varepsilon \in (0,\varepsilon^{\star})$ and any $w \in \mathcal{W}$, the closed-loop system possesses an isolated exponentially stable equilibrium point $(\bar{x},\bar{\eta})$ satisfying $\bar{e} = h(\bar{x},k(\bar{\eta}),w) = 0$. 
\end{theorem}

\medskip

Before proving the result, we examine how Theorem \ref{Thm:Desoer} reduces to known conditions in special cases. When restricted to weighted Euclidean vector norms $\|x\|_{\rm s} = (x^{\T}Px)^{1/2}$ where $P \succ 0$, \ref{Thm:Desoer-2} {\tb holds if and only if either of the following equivalent conditions hold \cite{MV:78,AP-AP-NVDW-HN:04}
\begin{subequations}
\begin{align}
\label{Eq:ContractionEuclidean}
\tfrac{\partial F_{\rm s}}{\partial \eta}(\eta,w)^{\T}P + P\tfrac{\partial F_{\rm s}}{\partial \eta}(\eta,w) \preceq -2\rho_{\rm s} P\\
\label{Eq:NewMonotone}
(\eta-\eta^{\prime})^{\T}P(F_{\rm s}(\eta,w)-F_{\rm s}(\eta^{\prime},w)) \leq -\rho_{\rm s}\|\eta-\eta^{\prime}\|_{P}^2
\end{align}
\end{subequations}
for some $\rho_{\rm s} > 0$, all $w \in \mathcal{W}$ and all $\eta, \eta^{\prime} \in \mathcal{C}_{w}$. If we restrict $k(\eta) = \eta$, $P = I_p$ and $\mathcal{C}_{w} = \real^p$, \eqref{Eq:ContractionEuclidean} and \eqref{Eq:NewMonotone} both reduce to the mapping $u \mapsto \pi(u,w)$ being \emph{strongly monotone} on $\real^p$, which is the Desoer--Lin condition \cite{CD-CL:85}.}

In the LTI case \eqref{Eq:LTI}--\eqref{Eq:LTIIntegral}, the equilibrium input-to-error mapping $\pi(u,w)$ is given explicitly by
\[
\pi(u,w) = G(0)u + G_{w}(0)w,
\]
where $G(0)$ and $G_w(0)$ are the DC gain matrices defined in \eqref{Eq:DCGains}. A linear integral feedback gain $u = K\eta$ obviously satisfies \ref{Ass:Controller-1}, and the dynamics \eqref{Eq:NonlinearSlow} reduce to
\[
\dot{\eta} = F_{\rm s}(\eta,w) = -G(0)K\eta - G_w(0)w.
\]
{\tb Considering again weighted} Euclidean norms, from \eqref{Eq:ContractionEuclidean} we see that \ref{Thm:Desoer-2} reduces to the existence of $P \succ 0$ and $\rho_{\rm s} > 0$ such that
\[
-(G(0)K)^{\T}P - P(G(0)K) \preceq -2\rho_{\rm s} P.
\]
By standard Lyapunov results, this holds if and only if $-G(0)K$ is Hurwitz, and we therefore properly recover the classical Davison/Morari result for LTI systems.


{\tr

}



\medskip


\begin{pfof}{Theorem \ref{Thm:Desoer}} {\tb The proof is divided into several steps.}

{\tb \emph{Step \#1 \textemdash{} Equilibrium Analysis and Error Dynamics:}} Let $w \in \mathcal{W}$. Closed-loop equilibria $(\bar{x},\bar{\eta})$ are characterized by the equations
\begin{equation}\label{Eq:ClosedLoopEq}
0 = f(\bar{x},\bar{u},w), \quad 0 = h(\bar{x},\bar{u},w), \quad \bar{u} = k(\bar{\eta}).
\end{equation}
Given any $\bar{u} \in \mathcal{U}_{w}$, by \ref{Ass:Plant-1} the first equation in \eqref{Eq:ClosedLoopEq} can be solved for $\bar{x} = \pi_{x}(\bar{u},w)$; together, \ref{Ass:Plant-1}/\ref{Ass:Plant-3} imply that $\bar{x}$ is isolated. Eliminating $\bar{x}$ and $\bar{u}$ from \eqref{Eq:ClosedLoopEq}, we obtain the error-zeroing equation $0 = \pi(k(\bar{\eta}),w)$. From \ref{Thm:Desoer-1}--\ref{Thm:Desoer-2}, the dynamics \eqref{Eq:NonlinearSlow} are infinitesimally contracting on a forward-invariant convex set $\mathcal{C}_{w}$; it follows from the main contraction stability theorem (see, e.g., \cite{EDS:10}) that \eqref{Eq:NonlinearSlow} possess a unique equilibrium point $\bar{\eta} \in \mathcal{C}_{w}$, and hence $0 = \pi(k(\bar{\eta}),w)$ is uniquely solvable on $\mathcal{C}_{w}$. By \ref{Thm:Desoer-1}, $\bar{\eta}$ further satisfies $k(\bar{\eta}) \in \mathcal{U}_{w}$, which justifies the initial application of \ref{Ass:Plant-1}. Thus, there exists a unique closed-loop equilibrium $(\bar{x},\bar{\eta}) \in \mathcal{X} \times \mathcal{C}_{w}$.

{\tb
Define the new state variable
\[
\xi \define x - \pi_{x}(k(\eta),w)
\]
and the new time variable $\tau \define \varepsilon t$. With this, the dynamics \eqref{Eq:Plant},\eqref{Eq:Integral} may be written in singularly perturbed form as
\begin{equation}\label{Eq:ClosedLoopTimeScale}
\begin{aligned}
\frac{\mathrm{d}\eta}{\mathrm{d}\tau} &= -h(\xi + \pi_{x}(k(\eta),w),k(\eta),w)\\
&\define -\tilde{h}(\xi,\eta)\\
\varepsilon \frac{\mathrm{d}\xi}{\mathrm{d}\tau} &= f(\xi + \pi_{x}(k(\eta),w),k(\eta),w) + \varepsilon\frac{\partial \pi_x}{\partial u}\frac{\partial k}{\partial \eta}\tilde{h}(\xi,\eta)\\
&\define \tilde{f}(\xi,\eta)
\end{aligned}
\end{equation}
where we have suppressed the arguments of $\tfrac{\partial \pi_x}{\partial u}$ and $\tfrac{\partial k}{\partial \eta}$. The equilibrium point of interest is now $(\xi,\eta) = (\vzeros[n],\bar{\eta})$.
}

\smallskip

{\tb \emph{Step \#2 \textemdash{} Bounding the Slow Dynamics:}} Let $V_{\rm s}(\eta) = \tfrac{1}{2}\|\eta-\bar{\eta}\|_{\rm s}^2$, and for later use note that by equivalence of norms, there exist constants $0 < c_{2\mathrm{s}} \leq c_{2\mathrm{s}}^{\prime}$ such that $c_{2\mathrm{s}} \|z\|_2 \leq \|z\|_{\rm s} \leq c_{2\mathrm{s}}^{\prime} \|z\|_2$ for any $z \in \real^p$. {\tb We compute the upper right-hand derivative (see, e.g., \cite[Section 3.4]{HKK:02}) of $V_{\rm s}$ along \eqref{Eq:ClosedLoopTimeScale} as
\begin{equation}\label{Eq:LieNonSmooth}
D^+V_{\rm s}(\eta) = \|\eta-\bar{\eta}\|_{\rm s} \cdot D^{+}\|\eta-\bar{\eta}\|_{\rm s}
\end{equation}
where
\[
D^{+}\|\eta-\bar{\eta}\|_{\rm s} = \limsup_{\alpha\to 0^+} \frac{\|\eta-\bar{\eta}-\alpha \tilde{h}(\xi,\eta)\|_{\rm s}-\|\eta-\bar{\eta}\|_{\rm s}}{\alpha}.
\]
Since by \eqref{Eq:DefofPi} and \eqref{Eq:NonlinearSlow} we have that
\[
\begin{aligned}
F_{\rm s}(\eta,w) &= -\pi(k(\eta),w) = -h(\pi_{x}(k(\eta),w),k(\eta),w)\\
&= -\tilde{h}(0,\eta),
\end{aligned}
\]
we may write
\[
\begin{aligned}
\eta - \bar{\eta} -\alpha \tilde{h}(\xi,\eta) = \eta - \bar{\eta} +\alpha F_{\rm s}(\eta,w) + \alpha \Delta \tilde{h}(\xi,\eta),
\end{aligned}
\]
where $\Delta \tilde{h}(\xi,\eta) = \tilde{h}(0,\eta) -  \tilde{h}(\xi,\eta)$. Inserting this into our expression for $D^{+}\|\eta-\bar{\eta}\|_{\rm s}$ and using the triangle inequality, we find that
\begin{equation}\label{Eq:Dini2}
\begin{aligned}
D^{+}\|\eta-\bar{\eta}\|_{\rm s} &\leq \limsup_{\alpha\to 0^+} \frac{\|\eta-\bar{\eta}+\alpha F_{\rm s}(\eta,w)\|_{\rm s}-\|\eta-\bar{\eta}\|_{\rm s}}{\alpha}\\
&\qquad + \|\Delta \tilde{h}(\xi,\eta)\|_{\rm s}.
\end{aligned}
\end{equation}
Let $J(\eta,w) = \tfrac{\partial F_{\rm s}}{\partial \eta}(\eta,w)$ denote the Jacobian matrix of $F_{\rm s}$. Since $\mathcal{C}_{w}$ is convex, $\bar{\eta} \in \mathcal{C}_{w}$, and $F_{\rm s}(\bar{\eta},w) = 0$, it follows from the multivariable mean value theorem (e.g., \cite[Theorem 6.21]{WR:76}) that
\begin{equation}\label{Eq:FMean}
F_{\rm s}(\eta,w) = \underbrace{\left[\int_{0}^{1}J(\bar{\eta} + \gamma (\eta-\bar{\eta}),w)\,\mathrm{d}\gamma\right]}_{\define J_{\rm avg}(\eta,w)}(\eta-\bar{\eta})
\end{equation}
for any $ \eta \in \mathcal{C}_{w}$. By \ref{Thm:Desoer-2} there exists $\rho_{\rm s} > 0$ such that $\mu_{\rm s}(J(\eta,w)) \leq -\rho_{\rm s}$ for all $\eta \in \mathcal{C}_{w}$, where $\mu_{\rm s}$ is the matrix measure associated with $\|\cdot\|_{\rm s}$. Since the matrix measure is a subadditive function \cite[Chap. 2.2.2]{MV:02}, it follows that
\begin{equation}\label{Eq:Javg}
\begin{aligned}
\mu_{\rm s}(J_{\rm avg}(\eta,w)) &\leq \int_{0}^{1} \mu_{\rm s}(J(\bar{\eta} + \gamma (\eta-\bar{\eta}),w))\,\mathrm{d}\gamma\\
& \leq -\rho_{\rm s} \int_{0}^{1}\mathrm{d}\gamma = -\rho_{\rm s}.
\end{aligned}
\end{equation}
%
%
Inserting \eqref{Eq:Dini2} and \eqref{Eq:FMean} into \eqref{Eq:LieNonSmooth} and using submultiplicativity of the induced matrix norm, we find that
\[
\begin{aligned}
D^{+}V_{\rm s}(\eta) &\leq \limsup_{\alpha\to 0^+} \frac{\|I_{p} + \alpha J_{\rm avg}(\eta,w)\|_{\rm s}-1}{\alpha}\|\eta-\bar{\eta}\|_{\rm s}^2\\
&\qquad + \|\Delta \tilde{h}(\xi,\eta)\|_{\rm s}\|\eta-\bar{\eta}\|_{\rm s}\\
&= -\mu_{\rm s}(J_{\rm avg}(\eta,w))\|\eta-\bar{\eta}\|_{\rm s}^2 + \|\Delta \tilde{h}(\xi,\eta)\|_{\rm s}\|\eta-\bar{\eta}\|_{\rm s}\\
&\leq -\rho_{\rm s}\|\eta-\bar{\eta}\|_{\rm s}^2 + \|\Delta \tilde{h}(\xi,\eta)\|_{\rm s}\|\eta-\bar{\eta}\|_{\rm s}.
\end{aligned}
\]
where we used \eqref{Eq:Javg} in the final inequality. Since $h$ is Lipschitz continuous, one easily finds that $\|\Delta \tilde{h}(\xi,\eta)\|_{\rm s} \leq c_{2\mathrm{s}}^{\prime}L_{h}\|\xi\|_{2}$ for some $L_h > 0$.
This yields the final bound
\begin{equation}\label{Eq:SlowBound}
D^{+}V_{\rm s}(\eta) \leq -\rho_{\rm s}\|\eta-\bar{\eta}\|_{\rm s}^2 + 2q_0 \|\xi\|_{2}\|\eta-\bar{\eta}\|_{\rm s},
\end{equation}
where $q_0 = \tfrac{1}{2}c_{2\mathrm{s}}^{\prime}L_{h}$.}

\smallskip

{\tb \emph{Step \#3 \textemdash{} Bounding the Fast Dynamics:} Begin by defining the deviation vector field $\map{g}{\real^n \times \mathcal{I}}{\real^n}$ by
\begin{equation}\label{Eq:gVectorField}
g(\xi,u,w) = f(\xi + \pi_x(u,w),u,w),
\end{equation}
For $r > 0$, let $B_{r}(0) = \setdef{\xi \in \real^n}{\|\xi\|_2 < r}$. Under Assumption \ref{Ass:Plant}, the converse Lyapunov result \cite[Lemma 9.8]{HKK:02} implies that there exist constants $r, c_1, c_2, \rho_{\rm f}, c_3, c_4 > 0$ and a continuously differentiable function $\map{V_{\rm f}}{B_{r}(0)\times\mathcal{I}}{\real}$ satisfying
\begin{subequations}\label{Eq:ConvLyap}
\begin{align}
\label{Eq:ConvLyap-1}
c_1 \|\xi\|_2^2 \leq V_{\rm f}(\xi,u,w) &\leq c_2 \|\xi\|_2^2\\
\label{Eq:ConvLyap-2}
\nabla_{\xi} V_{\rm f}(\xi,u,w)^{\T}g(\xi,u,w) &\leq -\rho_{\rm f} \|\xi\|_2^2\\
\label{Eq:ConvLyap-3}
\|\nabla_{\xi} V_{\rm f}(\xi,u,w)\|_2 &\leq c_3 \|\xi\|_2\\
\label{Eq:ConvLyap-4}
\|\nabla_{u} V_{\rm f}(\xi,u,w)\|_2 &\leq c_4 \|\xi\|_2^2.
\end{align}
\end{subequations}
for all $\xi \in B_{r}(0)$ and $(u,w) \in \mathcal{I}$. Along trajectories of \eqref{Eq:ClosedLoopTimeScale}, we compute that
\[
\begin{aligned}
\dot{V}_{\rm f}(\xi,k(\eta),w) &= T_1 + T_2\\
\end{aligned}
\]
where
\[
\begin{aligned}
T_1 &= \tfrac{1}{\varepsilon}\nabla_{\xi}V_{\rm f}(\xi,k(\eta),w)^{\T}\tilde{f}(\xi,\eta)\\
T_2 &= - \nabla_{u}V_{\rm u}(\xi,k(\eta),w)^{\T}\tfrac{\partial k}{\partial \eta}(\eta)\tilde{h}(\xi,\eta).
\end{aligned}
\]
Expanding $T_1$ using \eqref{Eq:ClosedLoopTimeScale} and \eqref{Eq:gVectorField}, we have
\[
\begin{aligned}
T_1 &= \tfrac{1}{\varepsilon}\nabla_{\xi}V_{\rm f}(\xi,k(\eta),w)^{\T}g(\xi,k(\eta),w)\\
&+ \nabla_{\xi}V_{\rm f}(\xi,k(\eta),w)^{\T}\frac{\partial \pi_x}{\partial u}\frac{\partial k}{\partial \eta}\tilde{h}(\xi,\eta).
\end{aligned}
\]
Using the equilibrium equation $\tilde{h}(0,\bar{\eta}) = 0$, we also have
\[
\begin{aligned}
\tilde{h}(\xi,\eta) &= \tilde{h}(\xi,\eta) - \tilde{h}(0,\bar{\eta}),
\end{aligned}
\]
and one quickly finds using \eqref{Eq:ClosedLoopTimeScale} that
\begin{equation}\label{Eq:htildeinequality}
\begin{aligned}
\|\tilde{h}(\xi,\eta)\|_2 &\leq L_h \left\|\begin{bmatrix}\xi+\pi_x(k(\eta),w)-\pi_x(k(\bar{\eta}),w) \\ k(\eta)-k(\bar{\eta})\end{bmatrix}\right\|_2\\
&\leq L_h \|\xi\|_2 + L_h L_k(L_{\pi_x}+1)\|\eta-\bar{\eta}\|_2
\end{aligned}
\end{equation}
where $L_{\pi_x}$ and $L_k$ are the Lipschitz constants of $\pi_x$ and $k(\eta)$, respectively, and where we used that $\|z\|_2 \leq \|z\|_{1}$ for $z \in \real^n$ to obtain the second inequality. We can now bound $T_1$ as
\[
\begin{aligned}
T_1 &\leq -\tfrac{1}{\varepsilon}\rho_{\rm f}\|\xi\|_2^2\\
&\quad  + c_3 L_{\pi_x} L_{k} \|\xi\|_2 (L_h \|\xi_k\| + L_h L_k(L_{\pi_x}+1)\|\eta-\bar{\eta}\|_2)
\end{aligned}
\]
where we used \eqref{Eq:ConvLyap-2} for the first term and \eqref{Eq:ConvLyap-3} in the second term. Similarly, we can use \eqref{Eq:htildeinequality} and \eqref{Eq:ConvLyap-4} to bound $T_2$ as
\[
\begin{aligned}
|T_2| &\leq c_4 L_k \|\xi\|_2^2(L_h \|\xi_k\|_2 + L_h L_k (L_{\pi_x}+1)\|\eta-\bar{\eta}\|_2)\\
&\leq c_4 L_kL_h r \|\xi\|_2^2 + c_4 L_h L_k^2 (L_{\pi_x}+1) r \|\xi_k\|_2\|\eta-\bar{\eta}\|_2
\end{aligned}
\]
where in the second line we used that $\xi \in B_{r}(0)$. Combining our inequalities, with minor manipulations we find that
\begin{equation}\label{Eq:FastBound}
\dot{V}_{\rm f}(\xi,k(\eta),w) \leq (\star)^{\T}\begin{bmatrix}
-\rho_{\rm f}/\varepsilon + q_1 & q_2\\ q_2 & 0
\end{bmatrix}
\begin{bmatrix}
\|\xi\|_2\\
\|\eta-\bar{\eta}\|_{\rm s}
\end{bmatrix},
\end{equation}
where the constants $q_1, q_2 > 0$ are independent of $\varepsilon$. Following \cite[Sec. 11.5]{HKK:02}, define the composite Lyapunov candidate
\[
V(\xi,\eta,w) = V_{\rm s}(\eta) + V_{\rm f}(x,k(\eta),w)
\]
for \eqref{Eq:ClosedLoopTimeScale}. Since $V_{\rm f}$ is continuous in all arguments, $k$ is continuous by \ref{Ass:Controller-1}, and $k(\bar{\eta}) \in \mathcal{U}_{w}$, it follows that $V_{\rm f}(\xi,k(\eta),w) > 0$ for all $(x,\eta) \in B_r(0) \times \mathcal{C}_{w}$ such that $(\xi,\eta) \neq (0,\bar{\eta})$. It follows immediately that $V$ is positive definite on $B_r(0) \times \mathcal{C}_{w}$ with respect to $(0,\bar{\eta})$. Combining the dissipation inequalities \eqref{Eq:SlowBound} and \eqref{Eq:FastBound}, we find that
\begin{equation}\label{Eq:Lie}
D^{+}V(\xi,\eta,w) \leq -\begin{bmatrix}
\|\xi\|_2 \\ 
\|\eta-\bar{\eta}\|_{\rm s}
\end{bmatrix}^{\T} Q \begin{bmatrix}
\|\xi\|_2 \\ 
\|\eta-\bar{\eta}\|_{\rm s}
\end{bmatrix}
\end{equation}
where
\[
Q = \begin{bmatrix}
\tfrac{\rho_{\rm f}}{\varepsilon}-q_1 & -(q_0 + q_2)\\
-(q_0 + q_2) & \rho_{\rm s}
\end{bmatrix}.
\]
Straightforward analysis shows that $Q \succ 0$ if and only if $\varepsilon \in (0,\varepsilon^{\star})$, where $\varepsilon^{\star} = \rho_{\rm f}\rho_{\rm s}/(\rho_{\rm s}q_1 + (q_0+q_2)^2) > 0$. Standard arguments using \eqref{Eq:ConvLyap-1} can now be applied to conclude that there exists a constant $\gamma > 0$ such that $D^{+}V(\xi,\eta,w) \leq -\gamma V(\xi,\eta,w)$ locally around $(0,\bar{\eta})$, and it quickly follows from the comparison lemma \cite[Lemma 3.4]{HKK:02} that the equilibrium is exponentially stable.
}
\end{pfof}

\smallskip


{\tb The conditions \ref{Thm:Desoer-1} and \ref{Thm:Desoer-2} in Theorem \ref{Thm:Desoer} ensure that the slow time-scale dynamics \eqref{Eq:NonlinearSlow} are infinitesimally contracting on a convex forward-invariant set $\mathcal{C}_{w}$, which by the contraction stability theorem (e.g., \cite{EDS:10}) yields the existence of a unique exponentially stable equilibrium point $\bar{\eta} \in \mathcal{C}_{w}$.\footnote{The additional minor assumption $k(\mathcal{C}_w) \subseteq \mathcal{U}_w$ in \ref{Thm:Desoer-1} ensures that the associated equilibrium control $\bar{u} = k(\bar{\eta})$ lies in the set of constant control inputs $\mathcal{U}_{w}$ from \eqref{Eq:Uw} for which there exists an equilibrium point for the state.} As is well-understood, a practical difficulty in applying contraction analysis comes in establishing the forward-invariance property. This can sometimes be done by exploiting structural properties of the vector field (e.g., for monotone systems). The situation simplifies considerably when the contraction condition \eqref{Eq:MatrixMeasureBound} can be \emph{globally} certified, in which case the forward-invariance conditions can be dropped. 

A situation of particular interest occurs when the equilibrium input-to-error mapping $\pi(u,w)$ has the separable form $\pi(u,w) = \pi_1(u) + \pi_2(w)$ for functions $\map{\pi_1}{\real^m}{\real^{p}}$ and $\map{\pi_2}{\real^{n_w}}{\real^p}$. For instance, $\pi$ will have this form when one considers reference tracking without exogenous disturbances. In this case, the reduced dynamics \eqref{Eq:NonlinearSlow} become
\begin{equation}\label{Eq:SeperableReducedDynamics}
\dot{\eta} = -\pi_1(k(\eta)) - \pi_2(w),
\end{equation}
which is similar to the LTI case in \eqref{Eq:SlowDynamics}. Mirroring the LTI case then, if $\pi_1$ is \emph{surjective} with right-inverse $\map{\pi_{1}^{\dagger}}{\real^p}{\real^m}$, then selecting $k = \pi_{1}^{\dagger}$ linearizes the dynamics \eqref{Eq:SeperableReducedDynamics} and ensures infinitesimal contraction. This provides a natural generalization of Davison's design $K = G(0)^{\dagger}$ to the nonlinear case. Note however that such an inversion-based design presupposes precise knowledge of $\pi_1$, along with the ability to compute a right-inverse. In Section \ref{Sec:SDP} we proceed down a different path for both global certification of stability and controller design, and present a LMI-based framework for synthesizing linear feedbacks $k(\eta) = K\eta$ when $\pi$ is only partially known.

\smallskip


\begin{remark}[\bf Alternatives to Contraction Conditions]\label{Rem:Alternatives}
The conditions \ref{Thm:Desoer-1}--\ref{Thm:Desoer-2} guarantee the existence/uniqueness of an equilibrium value for the integral state along with a stability property for the reduced dynamics \eqref{Eq:NonlinearSlow}. These assumptions can be separated and modified. For existence, one may assume that there exists a (continuously differentiable and Lipschitz continuous) equilibrium mapping $\map{c}{\mathcal{W}}{\real^p}$ such that $\pi(k(c(w)),w) = 0$ for all $w \in \mathcal{W}$. For stability, consider the deviation variable $\tilde{\eta} = \eta - c(w)$ and the corresponding dynamics
\begin{equation}\label{Eq:ReducedShifted}
\dot{\tilde{\eta}} = -\pi(k(\tilde{\eta}+c(w)),w),
\end{equation}
which now has an equilibrium at the origin for all $w \in \mathcal{W}$. The existence of $r, \mathsf{c}_1,\mathsf{c}_2,\mathsf{c_3},\mathsf{c}_4 > 0$ and a $w$-parameterized Lyapunov function $\map{V_{\rm s}}{B_{r}(0)\times\mathcal{W}}{\real}$ satisfying
\[
\begin{aligned}
\mathsf{c}_1\|\tilde{\eta}\|_2^2 \leq V_{\rm s}(\tilde{\eta},w) &\leq \mathsf{c}_2\|\tilde{\eta}\|_2^2\\
\nabla V_{\rm s}(\tilde{\eta},w)^{\T}[-\pi(k(\tilde{\eta}+c(w)),w)] &\leq -\mathsf{c}_3\|\tilde{\eta}\|_2^2\\
\|\nabla V_{\rm s}(\tilde{\eta},w)\|_2 &\leq \mathsf{c}_4 \|\tilde{\eta}\|_2
\end{aligned}
\]
establishes uniform (in $w$) exponential stability \cite[Sec. 9.6]{HKK:02} of the origin of \eqref{Eq:ReducedShifted}, and can be used in place of the contraction argument in Theorem \ref{Thm:Desoer}. As another option for a stability condition, one could instead build upon \eqref{Eq:NewMonotone}, and assume there exists a strongly convex and continuously differentiable function $\map{W}{\real^p}{\real_{\geq 0}}$ whose gradient is Lipschitz continuous on $\real^p$ and which satisfies
\begin{equation}\label{Eq:EIDMonotone}
\begin{aligned}
&(\nabla W(\eta) - \nabla W(\eta^{\prime}))^{\T}(\pi(k(\eta),w)-\pi(k(\eta^{\prime}),w))\\
&\qquad \qquad \qquad \qquad \quad \,\,\geq \rho_{\rm s} \|\eta-\eta^{\prime}\|_2^2
\end{aligned}
\end{equation}
for all $\eta,\eta^{\prime} \in \real^p$, all $w \in \mathcal{W}$, and some $\rho_{\rm s} > 0$. Inequalities of the form \eqref{Eq:EIDMonotone} arise, for example, when considering equilibrium-independent stability analysis \cite{JWSP:17e}. In this case, one may instead use $V_{\rm s}(\eta) = W(\eta) - W(\bar{\eta}) - \nabla W(\bar{\eta})^{\T}(\eta-\bar{\eta})$ in the proof of Theorem \ref{Thm:Desoer}, and the result goes through similarly. \hfill \oprocend
\end{remark}
}


\subsection{Application to Power System Frequency Control}
\label{Sec:PowerSystemFrequencyControl}

We illustrate the main result with a simple example arising in power system control. Our treatment is terse; we refer to \cite[Sec. 11.1.6]{PK:94} for engineering background and to \cite[Sec. IV]{FD-SB-JWSP-SG:19a},\cite{FD-SG:17,FD-SB-JWSP-SG:19a,JWSP:20b,JWSP:20d} for recent control-centric references. 

{\tb We consider an AC power system described by a model of the form \eqref{Eq:Plant} and satisfying Assumption \ref{Ass:Plant}; these assumptions are reasonable in practice, as low-level ``primary'' controllers in the system are designed to ensure stability \cite{FD-SB-JWSP-SG:19a}. The input $u \in \real^m$ represents changes to power injection set-points for controllable resources within the grid, while $w \in \real$ represents the net uncompensated load in the system. The measured error $e = \Delta \mathsf{f} \in \real$ is the deviation of the AC frequency from its nominal value (e.g., 60Hz). Many standard power system models have the property that the equilibrium input-to-error map \eqref{Eq:DefofPi} has the simple form
\begin{equation}\label{Eq:SteadyStateFrequency}
\Delta \mathsf{f} = \pi(u,w) = \tfrac{1}{\beta}\left(\sum_{i=1}^{m}\nolimits u_i - w\right)
\end{equation}
where $\beta > 0$ is called the \emph{frequency stiffness constant} of the system; see \cite[Sec. 11.1.6]{PK:94} for a derivation of \eqref{Eq:SteadyStateFrequency}.

The problem of \emph{secondary frequency regulation} is to design an integral control loop which regulates $\Delta\mathsf{f}$ to zero. A simple nonlinear integral control design to achieve this is
\begin{equation}\label{Eq:CAPI}
\tau \dot{\eta} = -\Delta \mathsf{f}, \quad u_i = \varphi_i(\eta),
\end{equation}
where $\eta \in \real$ is the integral state, $\tau > 0$ is the time constant, and $\map{\varphi_i}{\real}{\real}$ is a $L_i$-Lipschitz continuous function which satisfies the strong monotonicity condition
\begin{equation}\label{Eq:VarPhi}
(\varphi_i(\eta)-\varphi_i(\eta^{\prime}))(\eta-\eta^{\prime}) \geq \mu_i |\eta-\eta^{\prime}|^2
\end{equation}
for some $\mu_i > 0$ and all $\eta,\eta^{\prime} \in \real$. The nonlinearities $\varphi_i$ are used to optimally allocate the control resources; see, e.g., \cite{FD-SG:17} for details on this interpretation. With $\varepsilon = 1/\tau$, the controller \eqref{Eq:CAPI} is a special case of \eqref{Eq:Integral} with error signal $e = \Delta \mathsf{f}$ and feedback gain $k_i(\eta) = \varphi_i(\eta)$ satisfying \ref{Ass:Controller-1}. 

Combining \eqref{Eq:CAPI} and \eqref{Eq:SteadyStateFrequency}, it now follows that the reduced dynamics \eqref{Eq:NonlinearSlow} described in Theorem \ref{Thm:Desoer} are given by
\begin{equation}\label{Eq:CAPIPi}
\dot{\eta} = -\pi(k(\eta),w) = -\tfrac{1}{\beta}\sum_{i=1}^{m}\nolimits \varphi_i(\eta) + \tfrac{1}{\beta}w.
\end{equation}
For $\eta,\eta^{\prime} \in \real$ with $\pi = \pi(k(\eta),w)$ and $\pi^{\prime} = \pi(k(\eta^{\prime}),w)$, using \eqref{Eq:VarPhi} we have
\[
\begin{aligned}
(\eta-\eta^{\prime})(\pi -\pi^{\prime}) &= \tfrac{1}{\beta}\sum_{i=1}^{m}\nolimits(\varphi_i(\eta)-\varphi_i(\eta^{\prime}))(\eta-\eta^{\prime})\\
&\geq \left(\tfrac{1}{\beta}\sum_{i=1}^{m}\nolimits \mu_i\right)(\eta-\eta^{\prime})^2.
\end{aligned}
\]
We conclude that the contraction condition \eqref{Eq:NewMonotone} holds with $P = 1$ and $\rho_{\rm s} = \beta^{-1}\sum_{i}\nolimits \mu_i$. It follows that the slow dynamics \eqref{Eq:CAPIPi} are uniformly infinitesimally contracting on $\real$; \ref{Thm:Desoer-1}--\ref{Thm:Desoer-2} therefore both hold, and we conclude that the power system with controller \eqref{Eq:CAPI} is exponentially stable and achieves frequency regulation for sufficiently large $\tau > 0$. This extends the result of \cite{FD-SG:17} to general power system models, and allows for heterogeneity in the functions $\varphi_i$.
}

\section{Performance Analysis and Synthesis of Low-Gain Integral Controllers for Nonlinear and Uncertain Systems}
\label{Sec:SDP}

To complement the analytical results in Section \ref{Sec:Desoer}, we now pursue a computational framework for certifying performance of low-gain integral control schemes and synthesizing controller gains. To motivate our general approach, we return to the simple case of LTI systems in Section \ref{Sec:SDPLTI} before proceeding to nonlinear/uncertain system analysis and synthesis in Sections \ref{Sec:SDPNonlinear} and \ref{Sec:SDPNonlinearSynthesis}. Throughout we restrict our attention to linear feedbacks $k(\eta) = K\eta$. 


\subsection{Linear Time-Invariant Systems}
\label{Sec:SDPLTI}

We begin by considering a computational approach to the design of an integral feedback matrix $K$ in the LTI case; this will motivate our approach in subsequent sections.\footnote{The author is not aware of any literature applying LMI-based design techniques for low-gain output regulation, even for LTI systems. This is perhaps not surprising, given that the papers \cite{EJD:76,MM:85,CD-CL:85,HCT:92} preceded the development of computational methods for solving LMIs, and that output regulation research turned towards geometric methods in the wake of the seminal paper \cite{AI-CIB:90}.}
As described in Section \ref{Sec:Introduction}, the low integral feedback gain $\varepsilon$ induces a time-scale separation in the dynamics. In the LTI case, the slow dynamics are given by \eqref{Eq:SlowDynamics} with associated error output 
\begin{equation}\label{Eq:ErrorOutputLTI}
e(t) = G(0)K\eta(t) + G_w(0)w(t),
\end{equation}
and the sensitivity transfer matrix from $w$ to $e$ is therefore
\[
S_{\rm slow}(s) = s(sI_{p}+ \varepsilon G(0)K)^{-1}G_w(0).
\]
Davison's design recommendation \cite[Lemma 3]{EJD:76} $K_{\rm Dav} = G(0)^{\dagger}$ leads to the simple sensitivity function $S_{\rm slow}(s) = \frac{s}{s+\varepsilon}G_w(0)$, and achieves the minimum possible value
\[
\sup_{w \in \mathscr{L}_{2}[0,\infty), w \neq 0} \frac{\|e\|_{\mathscr{L}_2}}{\|w\|_{\mathscr{L}_2}} =  \|S_{\rm slow}(s)\|_{\mathcal{H}_{\infty}} = \|G_w(0)\|_2.
\]
%
%
for the induced $\mathscr{L}_2$-gain of the sensitivity function. This design however does not easily extend to the nonlinear case, may perform poorly in the presence of uncertainty (Figure \ref{Fig:StepResponseRobust}), and has the disadvantage for distributed linear control applications that $G(0)^{\dagger}$ is usually a dense matrix.

These disadvantages can be overcome by moving to a computational robust control framework. Due to the simple structure of the slow dynamics \eqref{Eq:SlowDynamics},\eqref{Eq:ErrorOutputLTI}, the design of $K$ can be formulated as an $\mathcal{H}_{\infty}$ state-feedback problem \cite[Chap. 7]{GED-FP:00}: find $Y \succ 0$ and $Z \in \real^{m \times p}$ such that
{\tb
\begin{equation}\label{Eq:LTISynthesis}
\begin{bmatrix}
    G(0)Z + (G(0)Z)^{\T} & \star & \star\\
    G_w(0)^{\T} & \gamma I_p  & \star\\
    -G(0)Z & -G_w(0) & \gamma I_p 
    \end{bmatrix} \succ 0,
\end{equation}
}
and then minimize over $\gamma > 0$. The resulting integral gain \textemdash{} which is recovered as $K_{\infty} = ZY^{-1}$ \textemdash{} will by construction achieve the same peak sensitivity as Davison's design, but the computational framework offers significant extensions. For instance, decentralization constraints $K \in \mathcal{K}$ where $\mathcal{K} \subseteq \real^{m \times p}$ is a subspace can be enforced by appending the additional constraints to \eqref{Eq:LTISynthesis} that $Y$ be diagonal and that $Z \in \mathcal{K}$.
%
%
%
%

\begin{figure}[ht!]
\centering
\begin{subfigure}{0.99\linewidth}
\includegraphics[width=\linewidth]{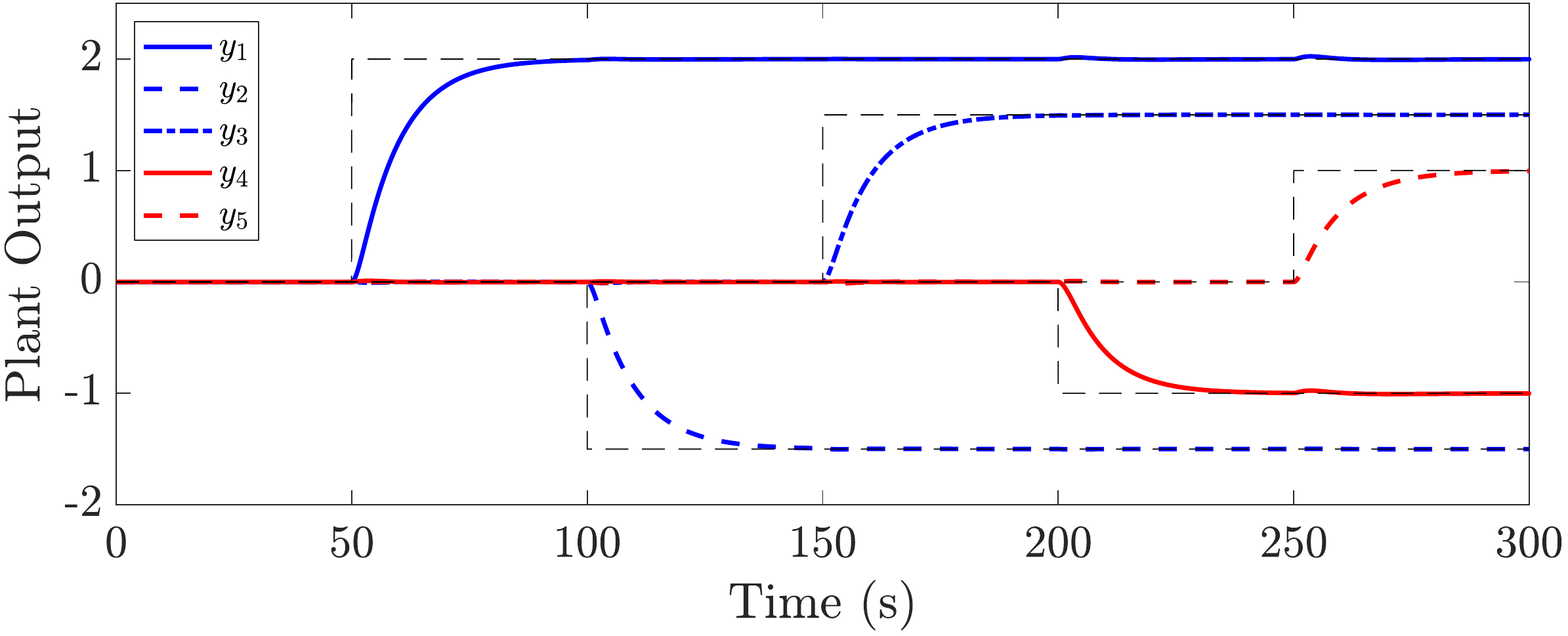}
\caption{Davison's selection $K = K_{\rm Dav} = G(0)^{\dagger}$, $\varepsilon = 0.1$.}
\label{Fig:StepResponse-a}
\end{subfigure}\\
\smallskip
\begin{subfigure}{0.99\linewidth}
\includegraphics[width=\linewidth]{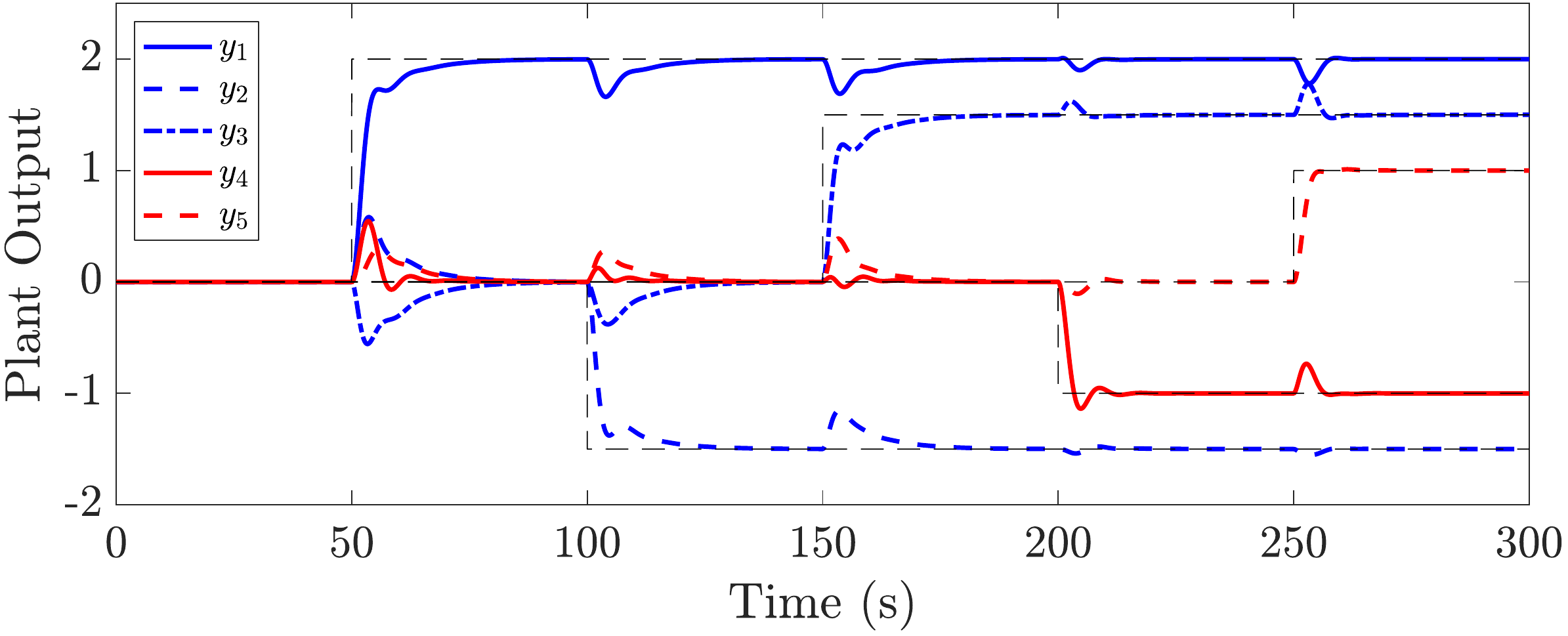}
\caption{Decentralized $\mathcal{H}_{\infty}$ design $K = K_{\infty}$, $\varepsilon = 3.5$.}
\label{Fig:StepResponse-b}
\end{subfigure}\\
\smallskip
\begin{subfigure}{0.99\linewidth}
\includegraphics[width=\linewidth]{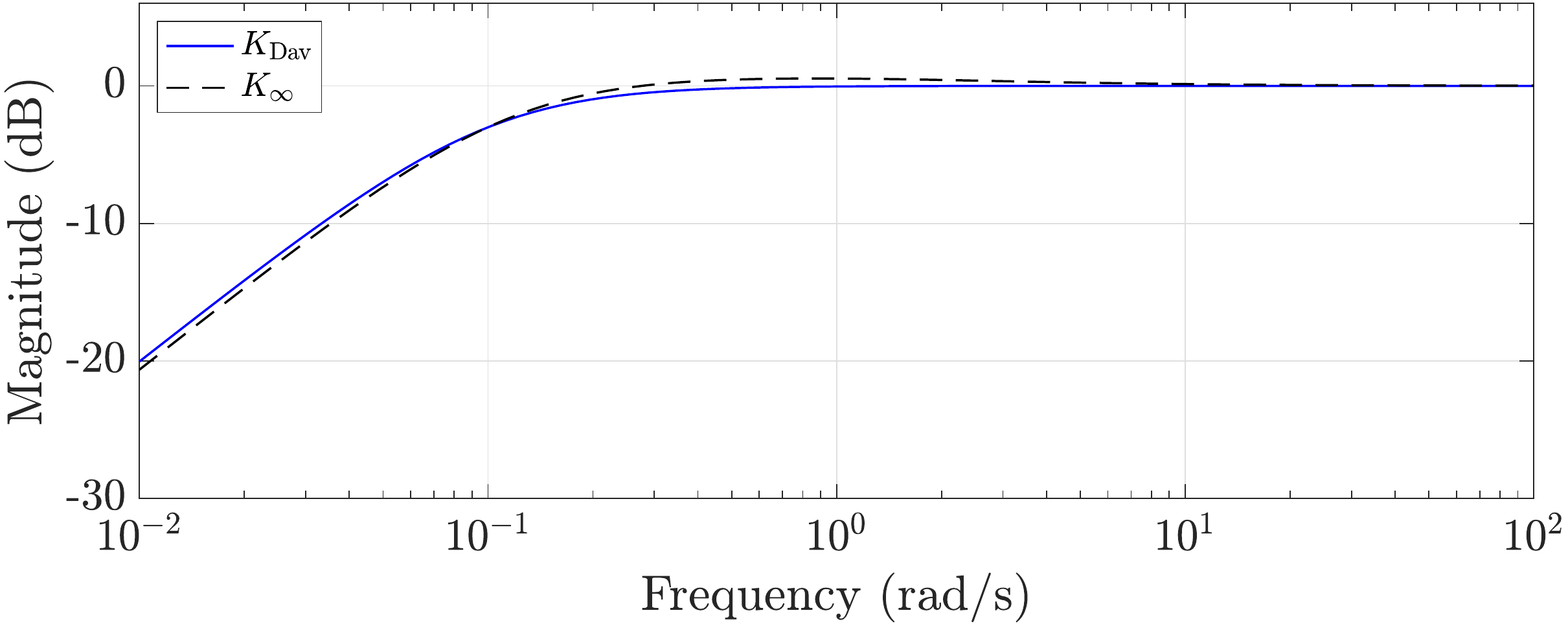}
\caption{Maximum singular value $\sigma_{\rm max}(S_{\rm slow}(\mathrm{j}\omega))$.}
\label{Fig:StepResponse-c}
\end{subfigure}
\caption{Reference step response of 7-input-5-output system with low-gain integral controllers.}
\label{Fig:StepResponse}
\end{figure}

We illustrate these ideas via reference tracking on a randomly generated stable LTI system with 30 states, 7 inputs, and 5 outputs.
%
We wish to design an feedback gain of the form 
\begin{equation}\label{Eq:Kstar}
K = \left[\begin{array}{@{}c|c@{}}
K_{11} & 0\\
\hline
0 & K_{22}
\end{array}\right], \quad K_{11} \in \real^{3\times 3}, \quad K_{22} \in \real^{4 \times 2},
\end{equation}
for use in a low-gain integral control scheme. The SDP \eqref{Eq:LTISynthesis} above was solved using SDPT3 with the YALMIP \cite{Lofberg2004} interface in MATLAB. Figures \ref{Fig:StepResponse-a},\ref{Fig:StepResponse-b} show the response of the resulting full-order closed-loop system to sequential step reference changes for the 5 output channels for the designs $K_{\rm Dav}$ and $K_{\infty}$, with associated maximum singular values of $S_{\rm slow}(\mathrm{j}\omega)$ plotted in Figure \ref{Fig:StepResponse-c}. The value of $\varepsilon$ was selected for the second design to match the bandwidth of the first design. The LMI-based design has no significant peaking in the sensitivity function and achieves the desired block-decentralization of the control actions.

\subsection{Nonlinear and Uncertain Systems}
\label{Sec:SDPNonlinear}


{\tb The slow time-scale reduced dynamics of the plant \eqref{Eq:Plant} and the controller \eqref{Eq:Integral} are as given in \eqref{Eq:NonlinearSlow}, which we rewrite here in input-output form as}
\begin{subequations}\label{Eq:SlowNonlinearDynamics}
\begin{align}
\label{Eq:SlowNonlinearDynamics-1}
\dot{\eta}(t) &= -\pi(K\eta(t),w(t)), \quad \eta(0) = \eta_0\\
\label{Eq:SlowNonlinearDynamics-2}
e(t) &= \pi(K\eta(t),w(t)).
\end{align}
\end{subequations}
{\tb The reduced dynamics \eqref{Eq:SlowNonlinearDynamics-1} can be expected to be significantly simpler and of lower dimension than the full nonlinear plant dynamics described by \eqref{Eq:Plant},\eqref{Eq:Integral}. As in Section \ref{Sec:SDPLTI}, our key observation is that the design of an optimal feedback gain $K$ in \eqref{Eq:SlowNonlinearDynamics} can be formulated as a state feedback control problem $u = K\eta$ for the integrator dynamics $\dot{\eta} = -\pi(u,w)$. This opens the door for applying analysis and computational techniques from robust control to synthesize integral gain matrices to achieve robust performance. To avoid some of the technical conditions  (e.g., forward-invariance) required by the general result in Section \ref{Sec:Desoer}, we will restrict our attention to global stability/performance analysis and synthesis for \eqref{Eq:SlowNonlinearDynamics}.}

{\tb Mirroring Section \ref{Sec:SDPLTI}, performance of \eqref{Eq:SlowNonlinearDynamics} will be quantified in terms of the energy contained in the disturbance and error signals.} To precisely formulate this, we let
\begin{equation}\label{Eq:SlowSignalSpace}
\map{\Sigma_{\eta_0}}{\mathscr{L}_{2e}^{n_w}[0,\infty)}{\mathscr{L}_{2e}^{p}[0,\infty)}, \quad e=\Sigma_{\eta_0}(w)
\end{equation}
denote the input-output signal-space operator defined by \eqref{Eq:SlowNonlinearDynamics}.\footnote{{\tb The subsequent LMI conditions we impose will be sufficient to guarantee that \eqref{Eq:SlowNonlinearDynamics} does indeed define such a mapping.}} We wish to establish performance guarantees which are independent of the operating point, which motivates the use of an \emph{incremental} $\mathscr{L}_2$-gain criteria \cite{CS-SW:15}. The system \eqref{Eq:SlowSignalSpace} is said to have incremental $\mathscr{L}_2$-gain less than or equal to $\gamma \geq 0$ if there exists $\beta \geq 0$ such that
\begin{equation}\label{Eq:IncrL2}
\begin{aligned}
\int_{0}^{T} \|e(t)-e^{\prime}(t)\|_2^2\,\mathrm{d}t &\leq \beta\|\eta_0-\eta_0^{\prime}\|_2^2\\
&\quad  + \gamma^2 \int_{0}^{T} \|w(t)-w^{\prime}(t)\|_2^2\,\mathrm{d}t
\end{aligned}
\end{equation}
for all $T > 0$, all initial conditions $\eta_0,\eta_0^{\prime}$, and all inputs $w,w^{\prime} \in \mathscr{L}_{2e}^{n_w}[0,\infty)$ where $e = \Sigma_{\eta_0}(w)$ and $e^{\prime} = \Sigma_{\eta_0^{\prime}}(w^{\prime})$. {\tb For LTI systems, this reduces to a standard $\mathcal{H}_{\infty}$ performance criterion.}

\begin{figure}[ht!]
\begin{center}
\tikzstyle{block} = [draw, fill=white, rectangle, 
    minimum height=3em, minimum width=6em, blur shadow={shadow blur steps=5}]
\tikzstyle{hold} = [draw, fill=white, rectangle, 
    minimum height=3em, minimum width=4em, blur shadow={shadow blur steps=5}]
\begin{tikzpicture}[auto, scale = 1, 
					node distance=2cm,
					>=latex',
					every node/.style={scale=1}]
    \node [hold] (system) {$\begin{array}{c|c|c}
    J & H & E_2 \\ \hline 
    G & F & E_1
    \end{array}$};
    
	\node [hold, above of=system, node distance=1.5cm] (delta) {$\Delta$};
    
    \coordinate [left of=system, node distance=1.7cm, yshift=0.24cm] (tmp-i-p);
    \coordinate [right of=system, node distance=1.7cm, yshift=0.24cm] (tmp-o-q);
    \draw [thick,-latex] (system.168) -| (tmp-i-p) |- node[pos=0.25, left] {$q$} (delta.west);
    \draw [thick,-latex] (delta.east) -| node[pos=0.75, right] {$p$} (tmp-o-q) |- (system.12);
    
    \coordinate [left of=system, node distance=2.2cm, yshift=-0.24cm] (tmp-o-z);
    \coordinate [right of=system, node distance=2.2cm, yshift=-0.24cm] (tmp-i-w);
    \draw [thick,-latex] (system.192) -- node[left, xshift=-0.5cm] {$\tilde{e} = \pi_{\Delta}(u,w)$} (tmp-o-z);
    \draw [thick,-latex] (tmp-i-w) node[right] {$\mathrm{col}(u,w)$}  -- (system.348);
    
    
    
\end{tikzpicture}
\end{center}
\caption{{\tb The linear fractional modelling framework.}}
\label{Fig:LFR}
\end{figure}
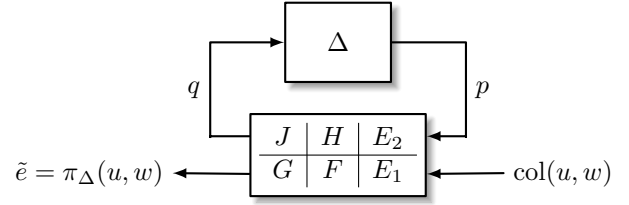

{\tb To work towards establishing \eqref{Eq:IncrL2}, we next represent the dynamics \eqref{Eq:SlowNonlinearDynamics} in standard \emph{linear fractional} form, i.e., in terms of a LTI system in feedback interconnection with an uncertain/nonlinear operator. We focus on representing the function $\pi(u,w)$ in this form, and} define a new function $\map{\pi_{\Delta}}{\real^m \times \real^{n_w}}{\real^p}$ via
\begin{equation}\label{Eq:LFR-1}
\begin{aligned}
\tilde{e} &= Fu + G\mathsf{p} + E_1w\\
\mathsf{q} &= Hu + J\mathsf{p} + E_2w\\
\mathsf{p} &= \Delta(\mathsf{q})
\end{aligned}
\end{equation}
where $(F,G,H,J,E_1,E_2)$ are fixed specified matrices, $\mathsf{p} \in \real^{n_{\mathsf{p}}}$ and $\mathsf{q} \in \real^{n_{\mathsf{q}}}$ are auxiliary variables, and $\map{\Delta}{\real^{n_{\mathsf{q}}}}{\real^{n_{\mathsf{p}}}}$ is a function. It is implicitly assumed that the representation \eqref{Eq:LFR-1} is \emph{well-posed}, in the sense that the mapping $\mathsf{q} \mapsto \mathsf{q} - J\Delta(\mathsf{q})$ is invertible on $\real^{n_{\mathsf{q}}}$; this holds trivially if $J = 0$. {\tb A diagram illustrating this functional representation is shown in Figure \ref{Fig:LFR}.}

We now allow $\Delta$ to range over a set $\boldsymbol{\Delta}$ of functions, which is coarsely described using point-wise incremental \emph{quadratic constraints}. In particular, we assume we have available a convex cone of symmetric matrices $\boldsymbol{\Theta} \subset \mathbb{S}^{(n_{\mathsf{p}} + n_{\mathsf{q}})\times(n_{\mathsf{p}} + n_{\mathsf{q}})}$ such that
\begin{equation}\label{Eq:IQC}
\begin{bmatrix}
\Delta(\mathsf{q}) - \Delta(\mathsf{q}^{\prime})\\
\mathsf{q} - \mathsf{q}^{\prime}
\end{bmatrix}^{\T}\Theta\begin{bmatrix}
\Delta(\mathsf{q}) - \Delta(\mathsf{q}^{\prime})\\
\mathsf{q} - \mathsf{q}^{\prime}
\end{bmatrix} \geq 0
\end{equation}
for all $\mathsf{q},\mathsf{q}^{\prime} \in \real^{n_{\mathsf{q}}}$, all $\Theta \in \boldsymbol{\Theta}$, and all $\Delta \in \boldsymbol{\Delta}$; the constraint $\Theta \in \boldsymbol{\Theta}$ must admit an LMI description. As a simple example, the cone of matrices
\[
\boldsymbol{\Theta} = \left\{\theta \left[\begin{smallmatrix}
-I_{n_{\mathsf{p}}} & 0\\
0 & L^2 I_{n_{\mathsf{q}}}
\end{smallmatrix}\right]\,\,\Big|\,\,\theta \geq 0\right\}
\]
can be used to describe nonlinearities $\Delta$ which are globally Lipschitz continuous with parameter $L$. We finally define
\begin{equation}\label{Eq:LFR}
\Pi \define \setdef{\pi_{\Delta}}{\Delta \in \boldsymbol{\Delta}}
\end{equation}
as the set of all maps one obtains. {\tb The modelling objective is to select $(F,G,H,J,E_1,E_2)$ and $\boldsymbol{\Delta}$ such that $\pi \in \Pi$, with as little conservatism introduced as possible.} While further explanation of this modelling framework is beyond our scope, we refer the reader to \cite{LEG-SN:00,CS:06,AB-LE-AN:09,CS-SW:15} for details, and we will illustrate with examples shortly. {\tb The LTI case \eqref{Eq:ErrorOutputLTI} is easily recovered by setting $F = G(0)$, $E_1 = G_w(0)$, and all other pieces of data equal to zero.}

%

{\tb We can now state the main analysis result, which certifies contraction and incremental $\mathscr{L}_2$-performance of the dynamics \eqref{Eq:SlowNonlinearDynamics} for all functions $\pi \in \Pi$.} For notational convenience, with $\gamma \geq 0$ we set
\[
\Theta_{\gamma} = \left[\begin{smallmatrix}
-\gamma^2 I_{n_w} & 0\\
0 & I_{p}
\end{smallmatrix}\right].
\]

\begin{theorem}[\bf Robust Stability and Performance of Reduced Dynamics]\label{Thm:RobustIncrementalSectorPerformance}
Consider the set $\Pi$ defined in \eqref{Eq:LFR-1}--\eqref{Eq:LFR} where $\Delta$ satisfies \eqref{Eq:IQC}, and assume that $\pi \in \Pi$. If there exists $K \in \real^{m \times p}$,  $P \succ 0$, and $\Theta \in \boldsymbol{\Theta}$ such that
\begin{equation}\label{Eq:RobustPerf}
{\small
\def\arraystretch{1.2}
\begin{aligned}
&(\star)^{\sf T}
\left[\begin{array}{@{}cc|c|c@{}}
0 & P & 0 & 0\\
P & 0 & 0 & 0\\
\hline
0 & 0 & \Theta & 0\\
\hline
0 & 0 & 0 & \Theta_{\gamma}
\end{array}\right]
\left[\begin{array}{@{}ccc@{}}
I_{p} & 0 & 0\\
-FK & -G & -E_1\\
\hline
0 & I_{n_{\mathsf{p}}} & 0\\
HK & J & E_2\\
\hline
0 & 0 & I_{n_w}\\
FK & G & E_1
\end{array}\right] \prec 0,
\end{aligned}
}
\end{equation}
then the following statements hold:
\begin{enumerate}
\item \label{Thm:RobustIncrementalSectorPerformance-1} the system \eqref{Eq:SlowNonlinearDynamics-1} is uniformly infinitesimally contracting with respect to norm $\|x\| = (x^{\T}Px)^{1/2}$ on $\real^p$;
%
\item \label{Thm:RobustIncrementalSectorPerformance-2} the incremental $\mathscr{L}_2$-gain of the reduced dynamics \eqref{Eq:SlowNonlinearDynamics} is less than or equal to $\gamma$.
\end{enumerate}
\end{theorem}

\smallskip

{\tb The modelling framework described in \eqref{Eq:LFR-1}--\eqref{Eq:LFR} is similar to the global linearization set-up in \cite[Chap. 4.3]{SB-LEG-EF-VB:94} and to the convergent system analysis in \cite[Chap. 5]{AP-NVDW-HN:05}. The main differences are that \eqref{Eq:IQC} allows for more general quadratic descriptions of the nonlinear/uncertain components than the typical small-gain uncertainty model, and that the LMI in Theorem \ref{Thm:RobustIncrementalSectorPerformance} additionally captures performance in response to disturbances. In the language of \cite{AP-NVDW-HN:05}, the LMI \eqref{Eq:RobustPerf} establishes both input-to-state convergence of the reduced dynamics \eqref{Eq:SlowNonlinearDynamics} and incremental $\mathscr{L}_2$ performance. Note however that we are imposing these conditions only on the reduced dynamics \eqref{Eq:SlowNonlinearDynamics}; we do not require the full closed-loop system \eqref{Eq:Plant},\eqref{Eq:Integral} to be contractive or input-to-state convergent.
}

\begin{pfof}{Theorem \ref{Thm:RobustIncrementalSectorPerformance}}
Fix any $\Delta \in \boldsymbol{\Delta}$ and let $\pi_{\Delta}$ be defined by \eqref{Eq:LFR-1}. Let $\xi,\xi^{\prime} \in \real^p$ and $w,w^{\prime} \in \real^{n_w}$ be arbitrary, set $u \define K\xi$ and $u^{\prime} \define K\xi^{\prime}$, and correspondingly define $(\tilde{e},\mathsf{q},\mathsf{p})$ and $(\tilde{e}^{\prime},\mathsf{q}^{\prime},\mathsf{p}^{\prime})$ via \eqref{Eq:LFR-1}. From \eqref{Eq:LFR-1} we find that
\begin{equation}\label{Eq:SubtractedEquations}
\begin{aligned}
\tilde{e} - \tilde{e}^{\prime} &= FK(\xi-\xi^{\prime}) + G(\mathsf{p}-\mathsf{p}^{\prime}) + E_1(w-w^{\prime})\\
\mathsf{q} - \mathsf{q}^{\prime} &= HK(\xi-\xi^{\prime}) + J(\mathsf{p}-\mathsf{p}^{\prime}) + E_2(w-w^{\prime})\\
\mathsf{p} - \mathsf{p}^{\prime} &= \Delta(\mathsf{q}) - \Delta(\mathsf{q}^{\prime}).
\end{aligned}
\end{equation}
{\tb By strict feasibility of \eqref{Eq:RobustPerf}, there exists $\rho_{\rm s} > 0$ such that 
\[
{\small
\def\arraystretch{1.2}
\begin{aligned}
&(\star)^{\sf T}
\left[\begin{array}{@{}cc|c|c@{}}
0 & P & 0 & 0\\
P & 0 & 0 & 0\\
\hline
0 & 0 & \Theta & 0\\
\hline
0 & 0 & 0 & \Theta_{\gamma}
\end{array}\right]
\left[\begin{array}{@{}ccc@{}}
I_{p} & 0 & 0\\
-FK & -G & -E_1\\
\hline
0 & I_{n_{\mathsf{p}}} & 0\\
HK & J & E_2\\
\hline
0 & 0 & I_{n_w}\\
FK & G & E_1
\end{array}\right] \preceq \left[\begin{smallmatrix}-2\rho_{\rm s} P & 0 & 0\\ 0 & 0 & 0 \\ 0 & 0 & 0\end{smallmatrix}\right].
\end{aligned}
}
\]
}
Left and right multiplying this by $\mathrm{col}(\xi-\xi^{\prime},\mathsf{p}-\mathsf{p}^{\prime},w-w^{\prime})$ and using \eqref{Eq:SubtractedEquations}, we obtain
{\tb
\begin{equation}\label{Eq:DissLMI}
\begin{aligned}
& \begin{bmatrix}
\xi - \xi^{\prime}\\
\tilde{e} - \tilde{e}^{\prime}
\end{bmatrix}^{\T}
\begin{bmatrix}
0 & -P\\
-P & 0
\end{bmatrix}
\begin{bmatrix}
\xi - \xi^{\prime}\\
\tilde{e} - \tilde{e}^{\prime}
\end{bmatrix} + \begin{bmatrix}
\mathsf{p}-\mathsf{p}^{\prime}\\
\mathsf{q} - \mathsf{q}^{\prime}
\end{bmatrix}^{\T}
\Theta
\begin{bmatrix}
\mathsf{p}-\mathsf{p}^{\prime}\\
\mathsf{q} - \mathsf{q}^{\prime}
\end{bmatrix}\\
&+ \|\tilde{e}-\tilde{e}^{\prime}\|_2^2 - \gamma^2\|w-w^{\prime}\|_2^2 \leq -2\rho_{\rm s} (\xi-\xi^{\prime})^{\T}P(\xi-\xi^{\prime}).
\end{aligned}
\end{equation}
}
To show statement \ref{Thm:RobustIncrementalSectorPerformance-1}, select $w^{\prime} = w$, and note using \eqref{Eq:IQC} that the second  and third terms in the above inequality are non-negative. Since $\tilde{e} = \pi_{\Delta}(K\xi,w)$, it follows that
{\tb
\[
(\pi_{\Delta}(K\xi,w)-\pi_{\Delta}(K\xi^{\prime},w))^{\T}P(\xi-\xi^{\prime}) \geq \rho_{\rm s} \|\xi-\xi^{\prime}\|_{P}^2
\]
}
which establishes the contraction condition \eqref{Eq:NewMonotone}. Since $\Delta \in \boldsymbol{\Delta}$ was arbitrary and $\pi \in \Pi$, statement \ref{Thm:RobustIncrementalSectorPerformance-1} holds. To show \ref{Thm:RobustIncrementalSectorPerformance-2} consider the extended input-output dynamics
\begin{equation}\label{Eq:OverboundedDynamics}
\begin{aligned}
\dot{\xi} &= -\tilde{e}\\
\tilde{e} &= \pi_{\Delta}(K\xi,w)\\
\xi(0) &= \xi_0
\end{aligned} \quad \begin{aligned}
\dot{\xi}^{\prime} &= -\tilde{e}^{\prime}\\
\tilde{e}^{\prime} &= \pi_{\Delta}(K\xi^{\prime},w^{\prime})\\
\xi^{\prime}(0) &= \xi_0^{\prime}
\end{aligned}
\end{equation}
and define $V(\xi,\xi^{\prime}) = \|\xi-\xi^{\prime}\|_{P}^2$. Differentiating along trajectories of \eqref{Eq:OverboundedDynamics} and inserting \eqref{Eq:DissLMI}, we obtain
\[
\begin{aligned}
\dot{V}(\xi,\xi^{\prime}) &= -2(\tilde{e}-\tilde{e}^{\prime})^{\T}P(\xi-\xi^{\prime})\\
&\leq -2\rho_{\rm s} V(\xi,\xi^{\prime}) - \|\tilde{e}-\tilde{e}^{\prime}\|_2^2 + \gamma^2 \|w-w^{\prime}\|_2^2.
\end{aligned}
\]
Integrating from $0$ to time $T > 0$ we have
\[
\begin{aligned}
V(\xi(T),\xi^{\prime}(T)) &- V(\xi_0,\xi^{\prime}_0) \\
&\leq - \int_{0}^{T} \|\tilde{e}-\tilde{e}^{\prime}\|_2^2 + \gamma^2 \|w-w^{\prime}\|_2^2\,\mathrm{d}t.
\end{aligned}
\]
Since $V(\xi(T),\xi^{\prime}(T)) \geq 0$, the performance inequality \eqref{Eq:IncrL2} now follows with $\beta = \lambda_{\rm max}(P)$. Since $\Delta \in \boldsymbol{\Delta}$ was arbitrary and $\pi \in \Pi$, statement \ref{Thm:RobustIncrementalSectorPerformance-2} holds.
\end{pfof}

\medskip

For a fixed $K$, the matrix inequality of Theorem \ref{Thm:RobustIncrementalSectorPerformance} is affine in $(P,\Theta,\gamma^2)$ and the best upper bound on $\gamma$ can be computed via semidefinite programming. While we have formulated Theorem \ref{Thm:RobustIncrementalSectorPerformance} for the nonlinear model \eqref{Eq:SlowNonlinearDynamics}, it can also be applied if the reduced dynamics are described by an \emph{uncertain linear} model; the details are standard and are omitted. 

The main user effort in applying Theorem \ref{Thm:RobustIncrementalSectorPerformance} is to appropriately select $F,G,H,J,E_1,E_2$ and $\boldsymbol{\Delta}$ such that $\pi \in \Pi$. As a simple illustration of the ideas, we continue our example from Section \ref{Sec:PowerSystemFrequencyControl}. In this case, we may directly model the mapping $\pi \circ k$ in \eqref{Eq:CAPIPi} by selecting $K = 1$, $F = 0$, $G = \beta^{-1}$, $E_1 = -\beta^{-1}$, $H = 1$, $J = 0$, and $E_2 = 0$, with $\mathsf{q} = \eta$ and $\mathsf{p} = \Delta(\mathsf{q}) = \sum_{i=1}^{m}\varphi_i(\mathsf{q})$. {\tb As each function $\varphi_i$ is $\mu_i$-strongly monotone and $L_i$-Lipschitz continuous, the nonlinear mapping $\Delta$ satisfies \eqref{Eq:IQC} with
\[
\boldsymbol{\Theta} = \left\{\theta \begin{bmatrix}
-2 & \mu + L\\
\mu + L & -2\mu L
\end{bmatrix}\,\,\Big|\,\,\theta \geq 0\right\},
\]
where $\mu = \sum_{i=1}^{m} \mu_i$ and $L = \sum_{i=1}^{m} L_i$. The performance LMI \eqref{Eq:RobustPerf} reduces to the $3\times 3$ LMI
\[
\left[\begin{array}{c|c|c}
2 \mu L \theta & \star & \star\\
\hline
P\beta^{-1} - \theta(\mu + L) & 2\theta - \beta^{-2} & \star \\ 
\hline
-P\beta^{-1} & \beta^{-2} & \gamma^2 - \beta^{-2} 
\end{array}\right] \succ 0.
\]
In this relatively simple case, one can analytically analyze this LMI, and one finds that it is feasible in $P, \theta > 0$ if and only if
\[
\gamma > \gamma^{\star} \define \frac{1}{\beta}\sqrt{\frac{(\kappa+1)^2}{4\kappa}}, \qquad \kappa = \frac{L}{\mu}.
\]
It follows that $\gamma^{\star}$ is the best certifiable upper bound on the incremental $\mathscr{L}_2$-gain of the slow time-scale dynamics \eqref{Eq:CAPIPi}. } 

\subsection{Synthesizing Feedback Gains for Robust Performance}
\label{Sec:SDPNonlinearSynthesis} 


Theorem \ref{Thm:RobustIncrementalSectorPerformance} can be exploited for direct \emph{convex} synthesis of controller gains for uncertain and nonlinear systems. The methodology here is inspired by procedures for robust state-feedback synthesis, and requires the the additional restrictions that the matrix $\Theta$ in \eqref{Eq:IQC} be nonsingular and that with the block partitionings
\begin{equation}\label{Eq:ThetaRestrictions}
\Theta = \begin{bmatrix}
\Theta_{11} & \Theta_{12}\\
\Theta_{12}^{\T} & \Theta_{22}
\end{bmatrix}, \quad \tilde{\Theta} \define \Theta^{-1} = \begin{bmatrix}
\tilde{\Theta}_{11} & \tilde{\Theta}_{12}\\
\tilde{\Theta}_{12}^{\T} & \tilde{\Theta}_{22}
\end{bmatrix}
\end{equation}
the sub-blocks satisfy $\Theta_{22} \succeq 0$ and $\tilde{\Theta}_{11} \preceq 0$.
Beginning with \eqref{Eq:RobustPerf}, define $Y = P^{-1} \succ 0$ and perform a congruence transformation on \eqref{Eq:RobustPerf} with $\mathrm{diag}(Y,I_{n_{\sf p}},I_{n_w})$. Multiplying through and setting $Z \define KY$ one obtains the following equivalent problem: find $Y \succ 0$, $Z \in \real^{m \times p}$ and $\Theta \in \boldsymbol{\Theta}$ such that
\begin{equation}\label{Eq:RobSyn}
{\small
\begin{aligned}
&(\star)^{\sf T}
\left[\begin{array}{@{}cc|c|cc@{}}
0 & I_p & 0 & 0 & 0\\
I_p & 0 & 0 & 0 & 0\\
\hline
0 & 0 & \Theta & 0 & 0 \\
\hline
0 & 0 & 0 & -I_{n_w} & 0\\
0 & 0 & 0 & 0 & \tfrac{1}{\gamma^2}I_{p}\\
\end{array}\right]
\def\arraystretch{1.2}
\left[\begin{array}{@{}ccc@{}}
I_{p} & 0 & 0\\
-FZ & -G & -E_1\\
\hline
0 & I_{n_{\sf p}} & 0\\
HZ & J & E_2\\
\hline
0 & 0 & I_{n_w}\\
FZ & G & E_1
\end{array}\right] \prec 0.
\end{aligned}
}
\end{equation}
Applying the \emph{Dualization Lemma} \cite[Cor. 4.10]{CS-SW:15}, \eqref{Eq:RobSyn} is equivalent to \eqref{Eq:RobSynDual} which is now affine in all decision variables. {\tb In further contrast to \cite[Chap. 5]{AP-NVDW-HN:05}, wherein a full-order output regulator design problem for the nonlinear plant \eqref{Eq:Plant} is treated, the low-gain design approach here is for the reduced dynamics \eqref{Eq:SlowNonlinearDynamics}, is a state-feedback design as opposed to an output feedback design, and our approach does not impose that the full closed-loop system \eqref{Eq:Plant},\eqref{Eq:Integral} be contractive/convergent.}

\begin{figure*}[t!]
\begin{equation}\label{Eq:RobSynDual}
(\star)^{\sf T}
\left[\begin{array}{@{}cc|cc|cc@{}}
0 & -I_p & 0 & 0 & 0 & 0\\
-I_p & 0 & 0 & 0 & 0 & 0\\
\hline
0 & 0 & -\tilde{\Theta}_{11} & -\tilde{\Theta}_{12} & 0 & 0\\
0 & 0 & -\tilde{\Theta}_{21} & -\tilde{\Theta}_{22} & 0 & 0 \\
\hline
0 & 0 & 0 & 0 & I_{n_w} & 0\\
0 & 0 & 0 & 0 & 0 & -\gamma^2 I_{p}\\
\end{array}\right]
\def\arraystretch{1.15}
\left[\begin{array}{@{}ccc@{}}
(FZ)^{\T} & -(HZ)^{\T} & -(FZ)^{\T}\\
I_{p} & 0 & 0\\
\hline
G^{\T} & -J^{\T} & -G^{\T}\\
0 & I_{n_{\sf q}} & 0\\
\hline
E_1^{\T} & -E_2^{\T} & -E_1^{\T}\\
0 & 0 & I_{p}
\end{array}\right] \prec 0.
\end{equation}
\end{figure*}

To illustrate these analysis and synthesis concepts, we return to the reference tracking example of Section \ref{Sec:SDPLTI}, and augment the previously generated system with three additional I/O channels $\mathsf{p}, \mathsf{q} \in \real^3$ with randomly selected coefficients in $G,H,J,E_2$. These channels are subject to the interconnection
\begin{equation}\label{Eq:DeltaExample}
\mathsf{p} = \Delta(\mathsf{q}) = \mathrm{col}(\mathsf{sat}(\mathsf{q}_1),\delta \mathsf{q}_2,\delta \mathsf{q}_3)
\end{equation}
where $\delta \in [-1,1]$ is an uncertain real parameter and $\mathsf{sat}$ denotes the standard saturation function. The associated equilibrium input-to-error map $\pi$ is now both nonlinear and uncertain. For $\Delta$ given in \eqref{Eq:DeltaExample}, the constraint \eqref{Eq:IQC} holds with
\begin{equation}\label{Eq:ThetaExample}
\Theta = \mathrm{daug}\left(\theta_{\varphi}\begin{bmatrix}
-2 & 1\\
1 & 0
\end{bmatrix},\begin{bmatrix}
-Q & S\\
S^{\T} & Q
\end{bmatrix}
\right)
\end{equation}
where $\theta_{\varphi} > 0$, $S = \left[\begin{smallmatrix}0 & s \\ -s & 0\end{smallmatrix}\right]$, $s \in \real$, and $Q = \left[\begin{smallmatrix}q_{11} & q_{12} \\ q_{12} & q_{22}\end{smallmatrix}\right] \succ 0$.

\begin{figure}[ht!]
\centering
\begin{subfigure}{0.99\linewidth}
\includegraphics[width=\linewidth]{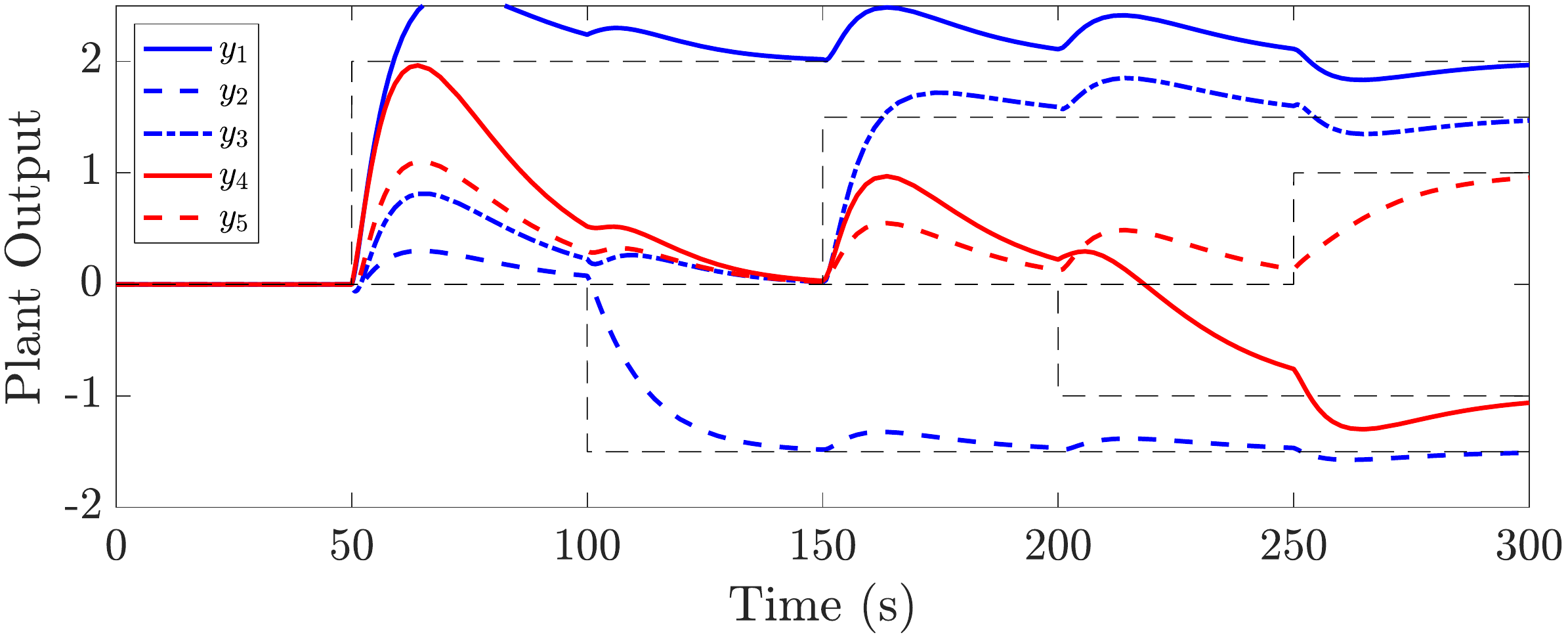}
\caption{Nominal design $K = F^{\dagger}$, $\varepsilon = 0.1$.}
\label{Fig:StepResponseRobust-a}
\end{subfigure}\\
\smallskip
\begin{subfigure}{0.99\linewidth}
\includegraphics[width=\linewidth]{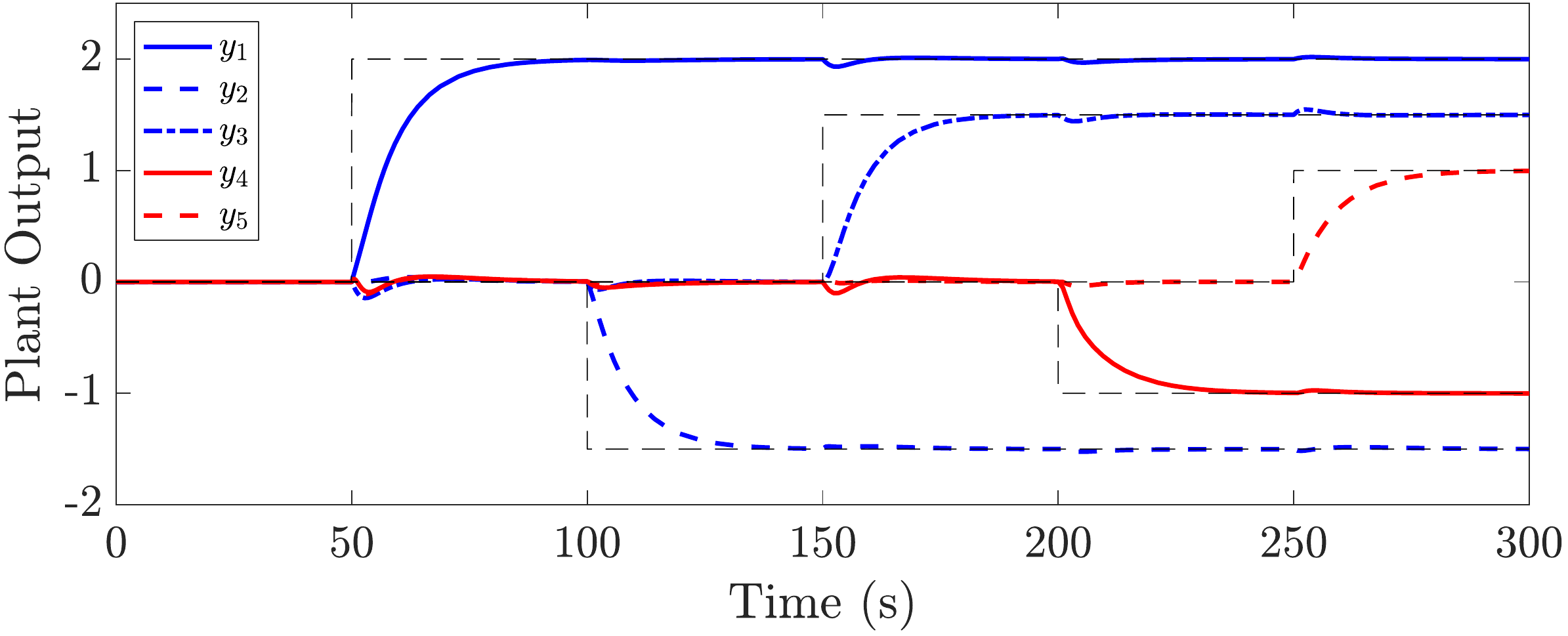}
\caption{Centralized robust design $K = K_{\rm robust}$, $\varepsilon = 3$.}
\label{Fig:StepResponseRobust-b}
\end{subfigure}\\
\smallskip
\begin{subfigure}{0.99\linewidth}
\includegraphics[width=\linewidth]{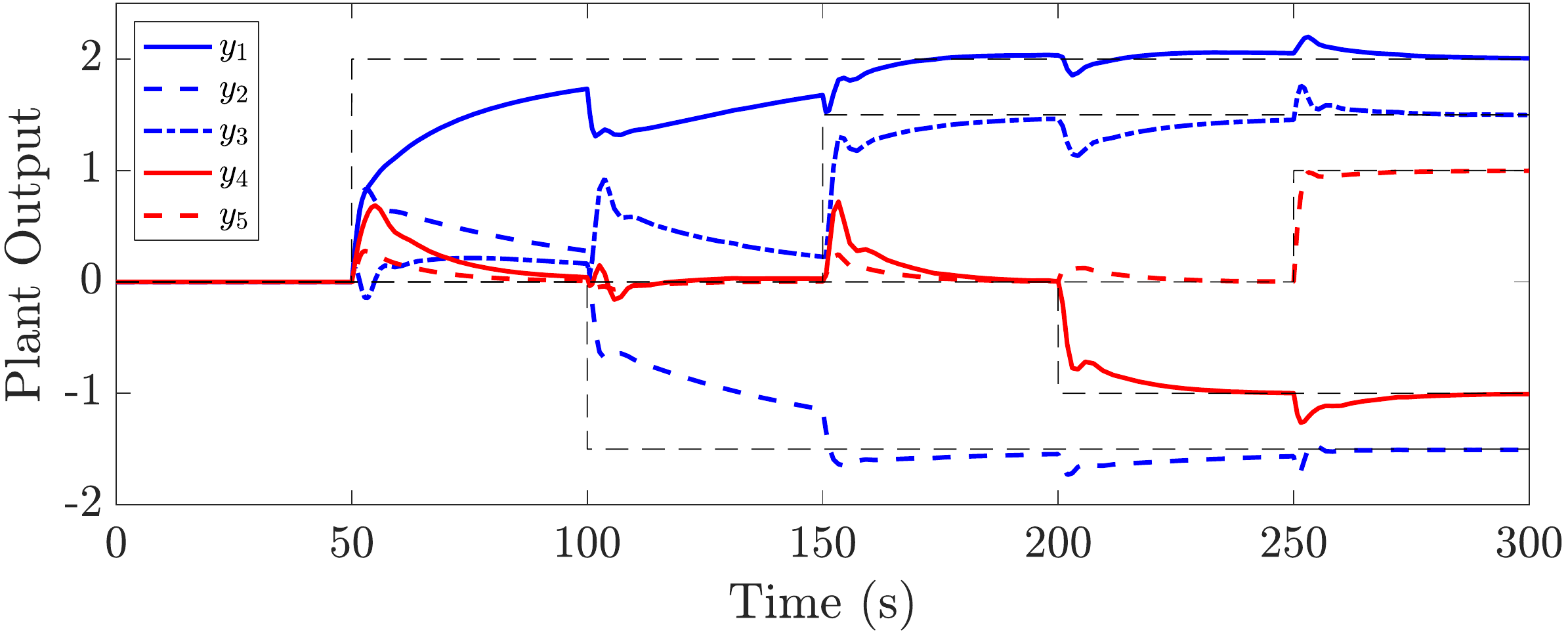}
\caption{Decentralized robust design $K = K_{\rm decent}$, $\varepsilon = 3$.}
\label{Fig:StepResponseRobust-c}
\end{subfigure}
\caption{Reference step response of nonlinear 7-input-5-output system with low-gain integral controllers.}
\label{Fig:StepResponseRobust}
\end{figure}

As a nominal controller design to compare against for this example, we use the same controller $K = G(0)^{\dagger} = F^{\dagger}$ from Section \ref{Sec:SDPLTI}, and attempt to certify robust performance via the SDP
\[
  \begin{aligned}
    \minimize_{\gamma,\, P \succ 0,\, \Theta \in \boldsymbol{\Theta}} &\quad \gamma^2 \quad \subto \quad \eqref{Eq:RobustPerf}.
  \end{aligned}
\]
Using SDPT3/YALMIP we can certify a performance bound of $\gamma = 16.9$ for the associated slow dynamics. While this is only an upper bound, it suggests that this nominal controller may perform poorly on instances of the nonlinear/uncertain system. Indeed,  Figure \ref{Fig:StepResponseRobust-a} shows the step response of the resulting full-order nonlinear closed-loop system to sequential step reference changes for the 5 output channels, for the case $\delta = -0.5$. The nominal design performs poorly in the presence of uncertainty and nonlinearity. To improve the design, we note that the matrix $\Theta$ in \eqref{Eq:ThetaExample} satisfies the additional restrictions mentioned in \eqref{Eq:ThetaRestrictions} if one restricts $s = 0$, in which case
\[
\tilde{\Theta} = \mathrm{daug}\left(\tilde{\theta}_{\varphi}\begin{bmatrix}
0 & 1\\
1 & 2
\end{bmatrix},\begin{bmatrix}
-D & 0\\
0 & D
\end{bmatrix}
\right)
\]
with $\tilde{\theta}_{\varphi} = \theta_{\varphi}^{-1} > 0$ and $D = \left[\begin{smallmatrix}d_{11} & d_{12} \\ d_{12} & d_{22}\end{smallmatrix}\right] = \left[\begin{smallmatrix}q_{11} & q_{12} \\ q_{12} & q_{22}\end{smallmatrix}\right]^{-1} \succ 0$. To synthesize a robust feedback gain, we solve the SDP
\[
  \begin{aligned}
    \minimize_{\gamma,\,Z,\, Y \succ 0,\, \tilde{\Theta} \in \boldsymbol{\tilde{\Theta}}} &\quad \gamma^2 \quad \subto \quad \eqref{Eq:RobSynDual},
  \end{aligned}
\]
and recover the controller as $K_{\rm robust} = ZY^{-1}$. For our example we obtain a much improved performance upper bound of $\gamma = 1.34$, and the step response, shown in Figure \ref{Fig:StepResponseRobust-b}, is significantly improved over the nominal design. If we wish to further impose a decentralization structure of the form \eqref{Eq:Kstar} on the controller, we solve the more constrained SDP
\[
  \begin{aligned}
    \minimize_{\gamma,\,Z \in \mathcal{K},\, Y \succ 0,\, \tilde{\Theta} \in \boldsymbol{\tilde{\Theta}}} &\quad \gamma^2 \quad \subto \quad \eqref{Eq:RobSynDual},\,\, Y\,\,\text{is diagonal},
  \end{aligned}
\]
{\tb and again recover the controller as $K_{\rm decent} = ZY^{-1}$.} For our example, a performance bound of $\gamma = 11.6$ is obtained, with step response shown in Figure \ref{Fig:StepResponseRobust-c}. The response is noticeably degraded over the centralized robust optimal design, but still improves on the nominal design, is robustly stable, and achieves decentralization of the integral action. 

As a final comment, robust analysis and synthesis procedures based on the (weighted) contraction norms $\|\cdot\|_{1}$ and $\|\cdot\|_{\infty}$ can also likely be developed, and may be valuable and efficient when the reduced dynamics \eqref{Eq:SlowNonlinearDynamics} is a monotone dynamical system; see \cite{SZK-CB-AR:15} for related ideas.

\section{Conclusions}
\label{Section: Conclusions}

Relaxed conditions have been given for stability of a nonlinear system under low-gain integral control, generalizing those available in the literature. The key idea is to impose an incremental-stability-type condition on the plant equilibrium input-to-error map. {\tb We then demonstrated how techniques from robust control can be applied to certify the key stability condition, and to synthesize integral controller gains using convex optimization which guarantee robust performance of the reduced dynamics. The results have been illustrated using analytical and numerical examples.}


{\tb Future work will focus on the application of these results to control problems in the energy systems domain and to the development of controllers for feedback-based optimization \cite{LSPL-ZEN-EM-JWSP:18e,LSPL-JWSP-EM:18l, MC-EDA-AB:18, SM-AH-SB-GH-FD:18, MC-JWSP-AB:19c}.} An open theoretical direction is to generalize the analysis and design approach presented here for tracking of signals generated by an arbitrary linear exosystem, which would yield a full generalization of \cite{EJD:76}. This generalization may require that incremental-type stability conditions be imposed on the plant, as considered in \cite{AP-LM:08}. {\tb Other open directions are to extend the low-gain results here to a discrete-time setup, and to anti-windup designs (e.g., \cite{GG-JH-IP-ART-MCT-LZ:03}), which should require only a modified analysis of the slow time-scale dynamics \eqref{Eq:NonlinearSlow}.}

\section*{Acknowledgements}
The authors thanks D. E. Miller, S. Jafarpour, F. Bullo, and A. Teel for helpful discussions.

%
%
%


\renewcommand{\baselinestretch}{1}
\bibliographystyle{IEEEtran}

\bibliography{/Users/jwsimpso/GoogleDrive/JohnSVN/bib/brevalias,%
/Users/jwsimpso/GoogleDrive/JohnSVN/bib/Main,%
/Users/jwsimpso/GoogleDrive/JohnSVN/bib/JWSP,%
/Users/jwsimpso/GoogleDrive/JohnSVN/bib/New%
}


\begin{IEEEbiography}[{\includegraphics[width=1in,height=1.25in,clip,keepaspectratio]{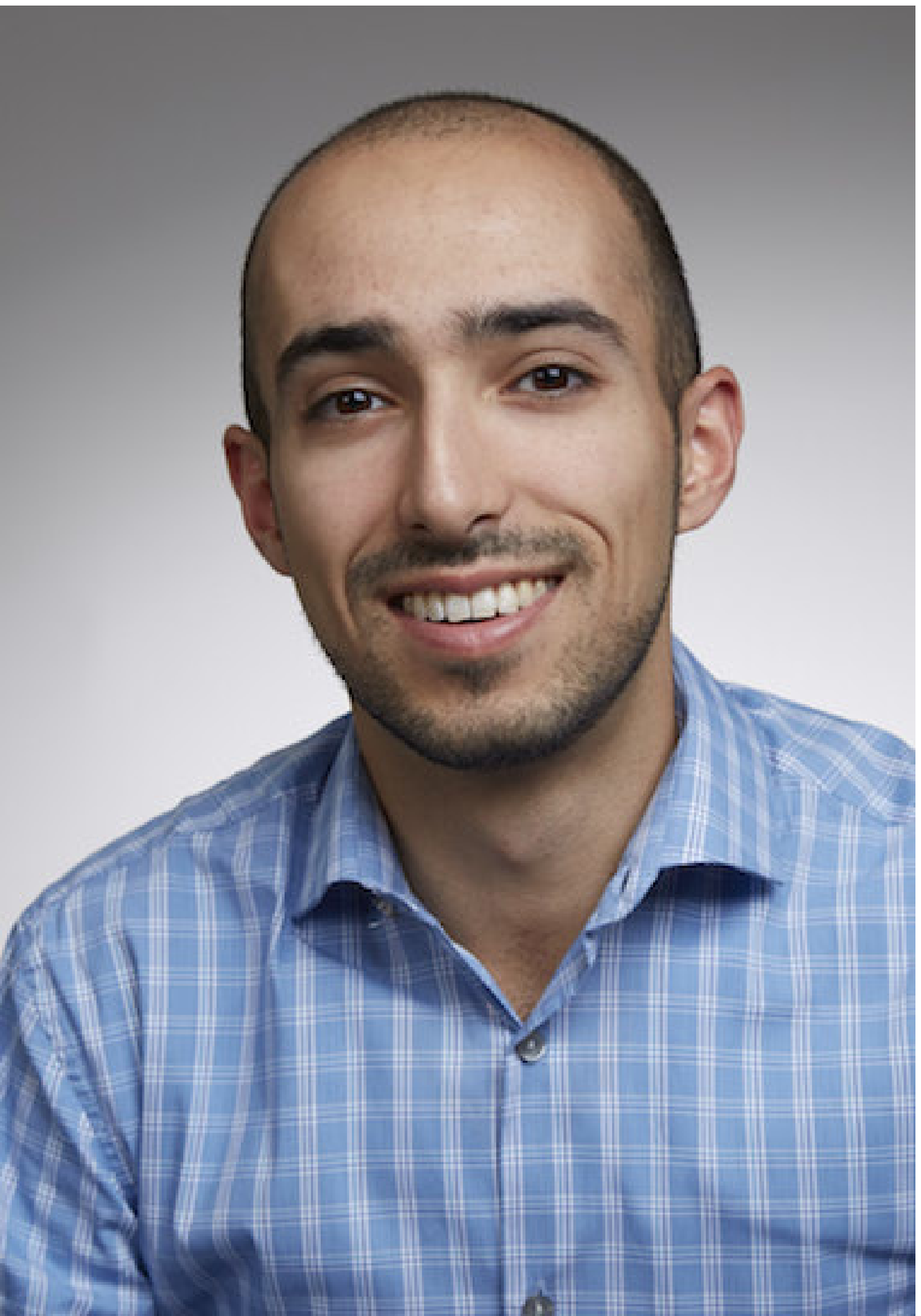}}]{John W. Simpson-Porco} (S'10--M'15--) received the B.Sc. degree in engineering physics from Queen's University, Kingston, ON, Canada in 2010, and the Ph.D. degree in mechanical engineering from the University of California at Santa Barbara, Santa Barbara, CA, USA in 2015.

He is currently an Assistant Professor of Electrical and Computer Engineering at the University of Toronto, Toronto, ON, Canada. He was previously an Assistant Professor at the University of Waterloo, Waterloo, ON, Canada and a visiting scientist with the Automatic Control Laboratory at ETH Z\"{u}rich, Z\"{u}rich, Switzerland. His research focuses on feedback control theory and applications of control in modernized power grids.

Prof. Simpson-Porco is a recipient of the Automatica Paper Prize, the Center for Control, Dynamical Systems and Computation Best Thesis Award, and the IEEE PES Technical Committee Working Group Recognition Award for Outstanding Technical Report. He is currently an Associate Editor for the IEEE Transactions on Smart Grid.
\end{IEEEbiography}

\end{document}